\documentclass[12pt]{amsart}

\usepackage{geometry}
\usepackage[utf8]{inputenc}
\usepackage[english]{babel}
\usepackage{amssymb,amsmath,amsthm}
\usepackage{enumitem}
\usepackage{amscd}

\theoremstyle{plain}
\newtheorem{theorem}{Theorem}[section]
\newtheorem{proposition}[theorem]{Proposition}
\newtheorem{lemma}[theorem]{Lemma}
\newtheorem{corollary}[theorem]{Corollary}
\newtheorem*{theorem*}{Theorem}

\theoremstyle{definition}
\newtheorem{definition}[theorem]{Definition}
\newtheorem{remark}[theorem]{Remark}
\newtheorem{example}[theorem]{Example}
\newtheorem{question}[theorem]{Question}

\newcommand{\N}{\mathbb{N}} 
\newcommand{\Z}{\mathbb{Z}} 
\newcommand{\R}{\mathbb{R}} 
\newcommand{\C}{\mathbb{C}} 
\newcommand{\D}{\mathbb{D}} 
\newcommand{\T}{\mathbb{T}} 
\newcommand{\K}{\mathbb{K}} 

\title{Zero-one law of orbital limit points for weighted shifts}
\author{Antonio Bonilla, Rodrigo Cardeccia, Karl-G. Grosse-Erdmann, and Santiago Muro}
\thanks{This publication is part of the project PID2022-139449NB-I00, funded by
MCIN/AEI/10.13039/501100011033/FEDER, UE
; the third author was supported by the Fonds de la Recherche Scientifique - FNRS under Grant n\textsuperscript{o} CDR J.0078.21; the second and fourth authors were supported by PICT 2018-4250 and CONICET}

\thanks{This work was partially done while the fourth author was visiting  K. Grosse-Erdmann at Universit\'e de Mons. He is extremely grateful for the hospitality and warmth with which he was received there.}

\address{Antonio Bonilla, Departamento de An\'alisis Matem\'atico and Instituto de Matem\'aticas y Aplicaciones (IMAULL),  Universidad de La Laguna, C/Astrof\'{\i}sico Francisco S\'anchez, s/n, 38721 La Laguna, Tenerife, Spain}
\email{abonilla@ull.edu.es}
\address{Rodrigo Cardeccia, Instituto Balseiro, Universidad Nacional de Cuyo – C.N.E.A. and CONICET, San Carlos de Bariloche,
Argentina}
\email{rodrigo.cardeccia@ib.edu.ar}
\address{Karl-G. Grosse-Erdmann, D\'epartement de Math\'ematique, Universit\'e de Mons, 20 Place du Parc, 7000 Mons, Belgium}
\email{kg.grosse-erdmann@umons.ac.be}
\address{Santiago Muro, FCEIA, Universidad Nacional de Rosario and CIFASIS, CONICET, Argentina}
\email{muro@cifasis-conicet.gov.ar}


\date{}

\begin{document}

\begin{abstract}
Chan and Seceleanu have shown that if a weighted shift operator on $\ell^p(\mathbb{Z})$, $1\leq p<\infty$, admits an orbit with a non-zero limit point then it is hypercyclic. We present a new proof of this result that allows to extend it to very general sequence spaces. In a similar vein we show that, in many sequence spaces, a weighted shift with a non-zero weakly sequentially recurrent vector has a dense set of such vectors; but an example on $c_0(\Z)$ shows that such an operator is not necessarily hypercyclic. On the other hand, we obtain that weakly sequentially hypercyclic weighted shifts are hypercyclic. Chan and Seceleanu have moreover shown that if an adjoint multiplication operator on a Bergman space admits an orbit with a non-zero limit point then it is hypercyclic. We extend this result to very general spaces of analytic functions, including the Hardy spaces.
\end{abstract}

\maketitle

\section{Introduction}

One of the most striking results in linear dynamics is due to Bourdon and Feldman \cite{BoFe03} who have shown that if an orbit of an operator is somewhere dense then it is already dense. Operators with a dense orbit are called hypercyclic. Clearly, if an orbit is only supposed to have a non-zero limit point one can no longer deduce that the orbit is dense. It therefore came as quite a surprise when Chan and Seceleanu \cite{Sec10}, \cite{ChSe12} showed that if a weighted shift operator has an orbit with a non-zero limit point, then the operator is necessarily hypercyclic.

\begin{theorem*}[Chan-Seceleanu zero-one law]\label{t-ChSeZeroOne}
Let $B_w$ be a unilateral \emph{(}respectively bilateral\emph{)}\linebreak weighted backward shift on $\ell^p$ \emph{(}respectively $\ell^p(\Z)$\emph{)}, $1\leq p<\infty$. 

If $B_w$ admits an orbit with a non-zero limit point then it is hypercyclic. 
\end{theorem*}

Let us emphasize however that an orbit with a non-zero limit point need not itself be dense; this aspect is further discussed in \cite{ChSe14}.

We note that the original proof of the zero-one law in the bilateral case is quite demanding. The starting point of our work was a new, shorter proof in that case. Also, Chan and Seceleanu do not discuss whether this behaviour extends to weighted shifts on other sequence spaces. Our proof allows us to obtain the zero-one law for all bilateral weighted shifts on all Fréchet sequence spaces in which the canonical sequence $(e_n)_n$ is a basis.

For the unilateral case, a corresponding extension has already been obtained by He, Huang and Yin \cite[Lemma 5]{HHY18}. A generalization of the theorem of Chan and Seceleanu to the notion of $\Gamma$-supercyclicity was recently provided by Abakumov and Abbar \cite{AbAb24}. 

In this paper we offer a thorough investigation of zero-one laws for unilateral and bilateral weighted backward shifts on general Fréchet sequence spaces. We follow Chan and Seceleanu \cite{ChSe12} in investigating whether the existence of limit points in a weaker sense already implies hypercyclicity. We obtain a number of positive results and some counter-examples, but we must leave some questions open for further research.

The paper is organized as follows. In Section \ref{s-ChanSecUni} we revisit the original zero-one law of Chan and Seceleanu for unilateral weighted backward shifts. Using essentially their proof we obtain a striking dichotomy: either the operator is hypercyclic, or else all orbits tend coordinatewise to zero (Theorem \ref{t-0-1-unilat_hyp}). 

Section \ref{s-ChanSecBil} is devoted to bilateral weighted backward shifts. We show there the zero-one law for such operators on all reasonable Fréchet sequence spaces, with an argument that is shorter than the one of Chan and Seceleanu (Theorem \ref{t-0-1-bilat}). If we weaken the topology then we still see that if the underlying space is $\ell^p(\Z)$, $1\leq p<\infty$, and if a weighted shift has an orbit with a non-zero weak sequential limit point, then it is hypercyclic (Corollary \ref{c-weakseqlp}). But an example shows that this result fails quite dramatically for the space $c_0(\Z)$ (Example \ref{ex-c0}). Still, on many spaces, including $c_0(\Z)$, the existence of a non-zero weakly sequentially recurrent vector implies that the weighted shift has a dense set of such vectors (Theorem \ref{t-weaksr}). Moreover, B\`es, Chan and Sanders \cite{BCS05} have shown that every weakly sequentially hypercyclic weighted shift on $\ell^p(\Z)$, $1\leq p<\infty$, or $c_0(\Z)$ is hypercyclic. This extends to all reasonable Fréchet sequence spaces (Corollary \ref{c-weakseqhc}).

In Section \ref{s-AdjMult} we look at a different class of operators. Chan and Seceleanu \cite{ChSe10} have shown that the zero-one law holds for the adjoints of multiplication operators on Bergman spaces. They implicitly asked if this is also the case when the operator is defined on a Hardy space. We give here a positive answer by proving zero-one laws for the adjoints of multiplication operators on many Banach spaces of analytic functions (Theorem \ref{t-01adjoint}). We also obtain a trichotomy for these operators on Bergman spaces and, under an additional condition, on Hardy spaces (Theorems \ref{teo: trichotomy bergman} and \ref{teo: trichotomy hardy}).

\subsection*{Notations and definitions}

A Fr\'echet sequence space over $\N_0$ is a Fr\'echet space that is a subspace of the space $\mathbb{K}^{\mathbb{N}_0}$ of all (real or complex) sequences and such that each coordinate functional $x=(x_n)_n\mapsto x_k$, $k\geq 0$, is continuous. Here, $\K$ stands for $\R$ or $\C$. The unit sequences are denoted by $e_n=(\delta_{n,k})_k$. A weight sequence is a sequence $w=(w_n)_{n\geq 1}$ of non-zero scalars. The unilateral weighted backward shift $B_w$ is then defined by $B_w(x_n)_{n\geq 0}=(w_{n+1}x_{n+1})_{n\geq 0}$; and $B$ denotes the shift with unit weights. Fr\'echet sequence spaces over $\Z$ and bilateral weighted backward shifts are defined analogously. It follows from the closed graph theorem that a weighted backward shift is continuous as soon as it maps a Fr\'echet sequence space into itself.

It will be convenient to write $[x]_n$ for the $n$-th coordinate of the sequence $x$; for example, $[B^nx]_p=x_{p+n}$. The support of a sequence $x$ is given by $\{n:x_n\neq 0\}$.

The topology of a Fréchet space can either be induced by an increasing separating sequence $(\|\cdot\|_r)_{r\geq 1}$ of seminorms, or by 
the corresponding F-norm $\|x\|=\sum_{r=1}^\infty \frac{1}{2^r}\max(1, \|x\|_r)$, see \cite{MeVo97}, \cite{KPR84}, \cite[Definition 2.9]{GrPe11}. Throughout this paper, $\|\cdot\|$ will denote an F-norm if $X$ is a Fréchet space, and a norm if $X$ is a Banach space. The open ball of radius $\rho$ around $x\in X$ is denoted by $D(x,\rho)$. 

\begin{example}\label{ex-spaces}
In the unilateral case, our main examples are the spaces $\ell^p=\ell^p(\N)$, $1\leq p<\infty$, of $p$-summable sequences, $c_0=c_0(\N)$ of null sequences, and $H(\C)$, the Fréchet space of entire functions, which can be regarded as a sequence space via the identification of an entire function with its sequence of Taylor coefficients.

Similarly, in the bilateral case, examples are the corresponding spaces $\ell^p(\Z)$, $1\leq p<\infty$, and $c_0(\Z)$. Less familiar will be the space $H(\C^*)$ of analytic functions on $\C^*=\C\setminus \{0\}$, which can be regarded as a sequence space by identifying $f(z)=\sum_{n\in\Z} a_n z^n$, $z\in\C^*$, with $(a_n)_{n\in\Z}$. This space is endowed with the topology of uniform convergence on compacta, induced by the seminorms $\|(a_n)_n\|_r = \sum_{n<0} |a_n|(1/r)^n+\sum_{n\geq 0} |a_n|r^n$, $r\geq 1$.
\end{example}

\subsection*{Linear dynamical properties}

Let $T$ be a (continuous, linear) operator on a Fr\'echet space $X$. 

A vector $x\in X$ is called \textit{hypercyclic} if its orbit is dense, that is, if for every non-empty open set $U\subset X$, $T^nx\in U$ for infinitely many $n\geq 0$. The operator $T$ is called hypercyclic if it admits a hypercyclic vector.

A vector $x\in X$ is called \textit{recurrent} if, for every neighbourhood $V$ of $x$, $T^nx\in V$ for infinitely many $n\geq 0$. The operator $T$ is called recurrent if it has a dense set of recurrent vectors.

A vector $y\in X$ is called a \textit{limit point} of the orbit of a vector $x\in X$ if, for every neighbourhood $V$ of $y$, $T^nx\in V$ for infinitely many $n\geq 0$. Note that definition and terminology of limit points vary in the literature. In our terminology, a periodic point is a limit point of the orbit.

If one considers $X$ under its weak topology, one has the analogous notions of \textit{weak hypercyclicity}, \textit{weak recurrence} and \textit{weak limit point}. 

We turn to the corresponding sequential notions. 

A vector $x\in X$ is called \textit{weakly sequentially hypercyclic} if, for every $y\in X$, there is a subsequence $(n_k)_k$ such that 
\[
T^{n_k}x\to y
\]
in the weak topology. The operator $T$ is called weakly sequentially hypercyclic if it admits a weakly sequentially hypercyclic vector.

A vector $x\in X$ is called \textit{weakly sequentially recurrent} if there is a subsequence $(n_k)_k$ such that 
\[
T^{n_k}x\to x
\]
in the weak topology. The operator $T$ is called weakly sequentially recurrent if it has a weakly sequentially dense set of weakly sequentially recurrent vectors. 

A vector $y\in X$ is called a \textit{weak sequential limit point} of the orbit of a vector $x\in X$ if there is a subsequence $(n_k)_k$ such that 
\[
T^{n_k}x\to y
\]
in the weak topology. 

If $X$ is a dual Banach space, then the corresponding weak* notions like \textit{weak*-hypercyclicity}, \textit{sequential weak*-hypercyclicity}, etc. are defined analogously.

\begin{remark}\label{rem-seqhc}
(a) In the metrizable topology of a Fréchet space $X$, the notions of hypercyclicity, recurrence and limit points are sequential. In particular, $y$ is a limit point of the orbit of $x$ if and only if there is a subsequence $(n_k)_k$ such that $T^{n_k}x\to y$ in $X$. 

(b) It is important to note that the weak (and weak*) sequential notions can also be understood in a weaker and possibly more appropriate way. For example, one might call a vector $x$ weakly sequentially hypercyclic for $T$ if the smallest weakly sequentially closed set containing the orbit of $x$ is the whole space; see \cite{Shk07}, \cite[Section 10.4]{BaMa09}. We prefer to work here with the more naive definitions above. 
\end{remark}

In the context of linear dynamics, recurrence was first studied in detail by Costakis et al. \cite{CoPa12}, \cite{CMP14}; see \cite{BGLP22}, \cite{CaMu22}, \cite{GrLo23}, \cite{CaMu22+} and \cite{GrLo22+} for recent developments; weak (sequential) recurrence was recently introduced by Amouch et al. \cite{ABBM23}. The consideration of limit points seems to be due to Chan and Seceleanu \cite{Sec10}, \cite{ChSe10}, \cite{ChSe12} and \cite{ChSe14}; see also \cite{Cha21} for a survey. 

For more on linear dynamics we refer to the textbooks \cite{BaMa09} and \cite{GrPe11}.

\subsection*{Conjugacies} Let us recall that every weighted backward shift operator $B_w$ is conjugate to the (unweighted) backward shift operator $B$ if one transfers the weight to the space; see \cite[Section 4.1]{GrPe11}.

More precisely, let $B_w$ be a unilateral weighted backward shift on a Fréchet sequence space $X$ over $\N_0$. We define the sequence $v=(v_n)_n$ by $v_n = \big(\prod_{\nu=1}^nw_\nu\big)^{-1}$, $n\geq 0$, and a new sequence space
\[
X_v = \{ (x_n)_{n} : (x_n v_n)_n\in X\}.
\]
Then the map $J_v: X_v \to X$, $(x_n)_n\mapsto (x_nv_n)_n$ is a vector space isomorphism, which induces in $X_v$ the structure of a Fréchet sequence space. Moreover, the following diagram commutes:
\[
\begin{CD}
X_v    @>B>>    X_v\\
@V{J_v}VV @VV{J_v}V \\
X    @>B_w>> X.
\end{CD}
\]
The same situation appears in the bilateral case. The sequence $v$ is then defined as above for nonnegative indices, while for $n<0$ one takes $v_n=\prod_{\nu=n+1}^0w_\nu$.

Thus $B_w:X\to X$ and $B:X_v\to X_v$ are conjugate operators, so that they have the same dynamics. We will therefore usually perform proofs for the unweighted shift $B$ and then transfer the result to $B_w$. This greatly simplifies the technicality of the arguments.

\section{Zero-one laws for unilateral weighted shifts}\label{s-ChanSecUni}

We begin by revisiting the Chan-Seceleanu zero-one law for the unilateral weighted backward shifts $B_w$ on $\ell^p$, $1\leq p<\infty$. It says that either no orbit has a non-zero limit point, or else $B_w$ is hypercyclic, meaning that one orbit has every vector in the space as limit point. 

In fact, using a similar idea as in the proof of \cite[Theorem 1.1]{ChSe12}, we obtain an interesting dichotomy, and that on almost arbitrary Fr\'echet sequence spaces.

\begin{theorem}\label{t-0-1-unilat_hyp}
Let $X$ be a Fr\'echet sequence space over $\N_0$ in which $(e_n)_n$ is a basis, and let $B_w$ be a weighted backward shift on $X$. Then the following dichotomy holds:
\begin{tabbing}
either  \=\kill
either  \> -- $B_w$ is hypercyclic,\\
or      \> -- for any $x\in X$, $B_w^nx\to 0$ in $\K^{\N_0}$.
\end{tabbing}
\end{theorem}

\begin{proof} By the above conjugacy argument we can assume that the operator is the unweighted backward shift $B$. Thus suppose that there is some $x\in X$ such that $B^nx\not\to 0$ in $\K^{\N_0}$. By the definition of convergence in $\K^{\N_0}$, there is some $p\geq 0$ such that $(x_{p+n})_n=([B^nx]_p)_n$ does not converge to 0; in other words, there is a subsequence $(n_k)_k$ and some $\delta>0$ such that, for all $k\geq 1$,
\[
|x_{p+n_k}|\geq \delta.
\]        
On the other hand, since $x\in X$ and $(e_n)_n$ is a basis, we have that
\[
x_{n}e_n\to 0
\]
in $X$ as $n\to\infty$. Altogether we find that
\[
e_{p+n_k} = \frac{1}{x_{p+n_k}}x_{p+n_k}e_{p+n_k}\to 0.
\]
Applying $p$ times the operator $B$ we see that $e_{n_k} \to 0$ in $X$, which implies that $B$ is hypercyclic, see \cite[Theorem 4.3]{GrPe11}.
\end{proof}

This result implies, in particular, the equivalence of conditions (i)-(v) in Chan and Seceleanu \cite[Theorem 1.1]{ChSe12}. More generally, we have the following.

\begin{corollary}\label{c-0-1-unilat_hyp}
Let $X$ be a Fr\'echet sequence space over $\N_0$ in which $(e_n)_n$ is a basis, and let $B_w$ be a weighted backward shift on $X$. Then the following assertions are equivalent:
\begin{enumerate}
\item[\rm (i)] $B_w$ is hypercyclic;
\item[\rm (ii)] $B_w$ has an orbit with a non-zero limit point in $X$;
\item[\rm (iii)] $B_w$ has an orbit with a non-zero weak sequential limit point in $X$;
\item[\rm (iv)] $B_w$ has an orbit with a non-zero weak limit point in $X$;
\item[\rm (v)] $B_w$ has an orbit with a non-zero limit point in $\K^{\N_0}$.
\end{enumerate}
\end{corollary}

If the basis is unconditional, the equivalence of (i) and (ii) is due to  He, Huang and Yin \cite[Lemma 5]{HHY18}.

The following consequences do not seem to have been stated before. Recall that if the basis $(e_n)_n$ is unconditional, then every weighted backward shift with a non-trivial periodic point is chaotic, but this is not necessarily true without unconditionality, see \cite[Theorem 8 and Example p.\ 65]{Gro00}.

\begin{corollary}\label{c-unilat_wh}
Let $X$ be a Fr\'echet sequence space over $\N_0$ in which $(e_n)_n$ is a basis. 

\emph{(a)} Every weakly hypercyclic weighted backward shift on $X$ is hypercyclic.

\emph{(b)} If a weighted backward shift on $X$ has a non-trivial periodic point then it is hypercyclic.
\end{corollary}

We can also extend an observation of Costakis, Manoussos and Parissis \cite[p.\ 1619]{CMP14} from $\ell^p$ to general Fréchet sequence spaces.

\begin{corollary}\label{c-0-1recunil}
Let $X$ be a Fr\'echet sequence space over $\N_0$ in which $(e_n)_n$ is a basis, and let $B_w$ be a weighted backward shift on $X$. Then the following assertions are equivalent:
\begin{enumerate}
\item[\rm (i)] $B_w$ is hypercyclic;
\item[\rm (ii)] $B_w$ is recurrent;
\item[\rm (iii)] $B_w$ has a non-zero recurrent vector;
\item[\rm (iv)] there exists a non-zero vector $y\in X$ so that, for every neighbourhood $V$ of $y$, there is some $x\in X$ such that $B_w^nx\in V$ for infinitely many $n\geq 0$.
\end{enumerate}
\end{corollary}

Note that, in condition (iv), one might demand that $x\in V$; thus we have here a local variant of the notion of \textit{almost recurrence} as introduced in \cite[Definition 3.2]{CaMu22+}.

\smallskip

For dual Banach sequence spaces one may consider weak*-limit points. Thus suppose that $X=Z^*$ is the dual of a Banach sequence space $Z$ over $\N_0$ under the natural pairing. If the canonical unit sequences $e_n$, $n\in \N_0$, belong to $Z$, then
$j_Ze_n=e_n^*$ is the $n$-th coordinate functional on $X$ (here $j_Z:Z\to Z^{**}=X^*$ is the canonical inclusion). Thus the coordinate functionals are weak*-continuous.

Suppose now that $B_w$ is not hypercyclic, and let $y$ be a weak*-limit point of the  orbit of $x$ under $B_w$. 
As a consequence of the weak*-continuity of the coordinate functionals, weak*-convergence implies coordinatewise convergence. Therefore we conclude, by Theorem \ref{t-0-1-unilat_hyp}, that $y=0$. We have shown the following.

\begin{corollary}\label{c-0-1-unilat*}
Let $X$ be a Banach sequence space over $\N_0$ in which $(e_n)_n$ is a basis.  Suppose also that, under the natural pairing, $X=Z^*$ for a Banach sequence space $Z$ which contains $(e_n)_n$, and let $B_w$ be a weighted backward shift on $X$. Then the following assertions are equivalent:
\begin{enumerate}
\item[\rm (i)] $B_w$ is hypercyclic;
\item[\rm (ii)] $B_w$ has an orbit with a non-zero weak*-sequential limit point in $X$;
\item[\rm (iii)] $B_w$ has an orbit with a non-zero weak*-limit point in $X$.
\end{enumerate}
\end{corollary}
Let us remark that the conditions in Corollary \ref{c-0-1-unilat*} are satisfied if $(e_n)_n$ is a boundedly complete basis in $X$, see \cite[Theorems 3.2.10, 3.2.12]{AlKa06}. 
Recall that a basis $(e_n)_n$ is boundedly complete if, for any $(x_n)_n\in\K^{\N_0}$, $\sup_{N\geq 0} \|\sum_{n=0}^N x_ne_n\|<\infty$ implies that $\sum_{n=0}^\infty x_ne_n$ converges in $X$.

\smallskip

There is one more equivalent condition in Chan and Seceleanu \cite[Theorem 1.1]{ChSe12}, namely their condition (vi). We can generalize their result to a large family of Banach sequence spaces; recall that a basis $(e_n)_n$ is called monotone if, for all $x\in X$, $\|x\| = \sup_{N\geq 0} \| \sum_{n=0}^N x_n e_n\|$.

\begin{theorem}\label{t-0-1-unilat_hyp2}
Let $X$ be a Banach sequence space over $\N_0$ in which $(e_n)_n$ is a monotone basis, and let $B_w$ be a weighted backward shift on $X$.  Then the following assertions are equivalent:
\begin{enumerate}
\item[\rm (i)] $B_w$ is hypercyclic;
\item[\rm (ii)] $B_w$ has an orbit with infinitely many members in an open ball whose closure avoids the origin; that is, there is a vector $x\in X$, a non-zero vector $y\in X$, and some $0<r<\|y\|$ such that $B_w^nx\in D(y,r)$ for infinitely many $n\geq 0$.
\end{enumerate}
\end{theorem}

\begin{proof}
Suppose that (ii) holds. We need to show that $B_w$ is hypercyclic. It follows from the property of the basis that there is some $N\geq 0$ such that
\begin{equation}\label{eq-biggerr}
\Big\| \sum_{n=0}^N y_n e_n\Big\|>r.
\end{equation}
On the other hand, by assumption, there is a subsequence $(n_k)_k$ such that, for all $k\geq 1$,
\[
\| B_w^{n_k}x-y \|<r.
\]
Again by the property of the basis, this implies that, for all $k\geq 1$,
\begin{equation}\label{eq-lessr}
\Big\| \sum_{n=0}^N ([B_w^{n_k}x]_n - y_n)e_n \Big\|<r.
\end{equation}

Now assume that $B_w$ is not hypercyclic. Then, by Theorem \ref{t-0-1-unilat_hyp}, we have that
\[
B_w^n x\to 0\ \text{in } \K^{\N_0}.
\]
By \eqref{eq-lessr}, this implies that
\[
\Big\| \sum_{n=0}^N y_ne_n \Big\|\leq r,
\]
which contradicts \eqref{eq-biggerr}. Hence $B_w$ is hypercyclic.
\end{proof}

In a Banach sequence space with basis $(e_n)_n$, the norm
\[
|\!|\!| x |\!|\!| = \sup_{N\geq 0} \Big\| \sum_{n=0}^N x_n e_n\Big\|
\]
is an equivalent norm for which $(e_n)_n$ is monotone, see \cite[Remark 1.1.6]{AlKa06}. Thus we have the following.

\begin{corollary}\label{c-0-1-unilat_hyp2}
Let $X$ be a Banach sequence space over $\N_0$ in which $(e_n)_n$ is a basis, and let $B_w$ be a weighted backward shift on $X$. Then the following assertions are equivalent:
\begin{enumerate}
\item[\rm (i)] $B_w$ is hypercyclic;
\item[\rm (ii)] there is a vector $x\in X$, a non-zero vector $y\in X$, and some $0<r<|\!|\!| y |\!|\!|$ such that $B_w^nx\in D_{|\!|\!| \cdot |\!|\!|}(y,r)$ for infinitely many $n\geq 0$.
\end{enumerate}
\end{corollary}

\section{Zero-one laws for bilateral weighted shifts}\label{s-ChanSecBil}

It was already noted by Chan and Seceleanu \cite{ChSe12} that there is no general zero-one law for weak limit points in the bilateral case. Indeed, Chan and Sanders \cite{ChSa04} had previously shown that any space $\ell^p(\Z)$, $2\leq p <\infty$, supports a weighted bilateral shift that is a weakly hypercyclic but not hypercyclic, see also \cite[Corollary 10.28]{BaMa09}. Interestingly, Shkarin \cite{Shk07} has shown that this phenomenon cannot arise if $1\leq p <2$, see also \cite[Proposition 10.20]{BaMa09} and Theorem \ref{t-0-1 bil ell_p weak p<2} below.

Thus, Corollary \ref{c-0-1-unilat_hyp} does not hold in the bilateral setting. We will see that versions of it are true in particular spaces or under additional assumptions. We will discuss each of the equivalent conditions from Corollary \ref{c-0-1-unilat_hyp} in the following subsections.

A characterization of hypercyclic bilateral shifts is well known, see \cite[Theorem 4.13]{GrPe11}. The following expresses the characterizing condition for the unweighted shift in a very handy way.

\begin{theorem}\label{t-hcbil}
Let $X$ be a Fr\'echet sequence space over $\Z$ in which $(e_n)_n$ is a basis, and let the unweighted backward shift $B$ be an operator on $X$. 

Then $B$ is hypercyclic if and only if, for any $p\geq 0$ and any $\varepsilon>0$, there is some $n\geq 0$ such that
\[
\|e_{p+n}\|<\varepsilon\ \ \text{and}\text\ \ \|e_{p-n}\|<\varepsilon.
\]
\end{theorem}

\begin{proof} 
The condition is necessary by \cite[Theorem 4.12]{GrPe11}. On the other hand, suppose that the condition holds. Making $\varepsilon$ smaller, if necessary, and applying the continuous operator $B$ up to $2p$ times, we obtain that, for any $p\geq 0$ and any $\varepsilon>0$, there is some $n\geq 0$ such that, for any $j=-p,\ldots,p$,
\[
\|e_{j+n}\| <\varepsilon \quad\text{and} \quad \|e_{j-n}\| <\varepsilon.
\] 
By considering $\varepsilon$ even smaller, one necessarily has that $n\geq p$. Then this condition implies that $B$ is hypercyclic, see \cite[Theorem 4.12]{GrPe11} and its proof.
\end{proof}

\subsection{Existence of limit points}\label{s-lim} 
Chan and Seceleanu have obtained their zero-one law in the spaces $\ell^p(\Z)$, $1\leq p <\infty$. Their result can be extended, and even strengthened, to arbitrary Fréchet sequence spaces under the usual assumption. The proof is surprisingly short.

\begin{theorem}\label{t-0-1-bilat}
Let $X$ be a Fr\'echet sequence space over $\Z$ in which $(e_n)_n$ is a basis, and let $B_w$ be a weighted backward shift on $X$. Then the following assertions are equivalent:
\begin{enumerate}
\item[\rm (i)] $B_w$ is hypercyclic;
\item[\rm (ii)] $B_w$ has an orbit with a non-zero limit point;
\item[\rm (iii)] there exists a non-zero vector $y\in X$ such that, for every neighbourhood $V$ of $y$, there is some $x\in X$ such that $B_w^nx\in V$ for infinitely many $n\geq 0$.
\end{enumerate}
\end{theorem}

\begin{proof}
It suffices to show that (iii) implies (i). As before, we can assume that the shift is unweighted.

Let $y\in X$ be a non-zero vector  satisfying the condition in (iii). In order to show that $B$ is hypercyclic, fix $p\geq 0$ and $\varepsilon>0$. Since $y\neq 0$, there is some $q\in \Z$ such that $y_q\neq 0$, and we can assume that $y_q=2$. 

Since $(e_n)_n$ is a basis, the projections $z\mapsto z_ne_n$, $n\in\Z$, are equicontinuous on $X$ by the Banach-Steinhaus theorem. There is therefore some $\eta>0$ such that, for any $z\in X$ with $\|z\|<\eta$ and any $n\in\Z$, 
\begin{equation}\label{eq-cond1}
\|z_ne_n\|< \varepsilon.    
\end{equation}
By making $\eta$ smaller, if necessary, we can assume that for any $z\in X,$
\begin{equation}\label{eq-cond2}
\|z-y\|< \eta \ \Longrightarrow \ |z_q|\geq 1.
\end{equation} 

By hypothesis there is then some $x\in X$ and a subsequence $(n_k)_k$ such that, for all $k\geq 1$,
\[
\|B^{n_k}x-y\|<\frac{\eta}{2},
\]
where we may assume that $n_1\geq 2p-2q$. By \eqref{eq-cond2}, it follows that, for all $k\geq 1$,
\begin{equation}\label{eq-cond3}
|x_{q+n_k}|=|[B^{n_k}x]_q|\geq  1.
\end{equation} 

By the basis assumption, there is also some $N\in\N$ such that $w:=\sum_{n=-N}^N y_ne_n$ satisfies $\|y-w\|<\frac{\eta}{2}$ and hence that, for all $k\geq 1$,
\begin{equation}\label{eq-cond4}
\|B^{n_k}x-w\|<\eta.
\end{equation}

Note that
\[
[B^{n_k}x-w]_{q+n_1-n_k}=x_{q+n_1} - w_{q+n_1-n_k}= x_{q+n_1}
\]
whenever $k$ is so large that $n_k>q+n_1+N$. Then \eqref{eq-cond1}, \eqref{eq-cond3} and \eqref{eq-cond4} imply that, for such $k$,
\[
\|e_{q+n_1-n_k}\|=\Big\|\frac{1}{x_{q+n_1}}[B^{n_k}x-w]_{q+n_1-n_k}e_{q+n_1-n_k}\Big\|<\varepsilon,
\]
where we have used that $\|az\|\leq \|z\|$ for any scalar $a$ with $|a|\leq 1$ and any $z\in X$, an inequality that follows from the definition of F-norm.

When we set
\[
n:= p-q+n_k-n_1,
\]
then we have for all sufficiently large $k$ that $n\geq 0$ and
\begin{equation}\label{eq-app1}
\|e_{p-n}\|=\|e_{q+n_1-n_k}\|<\varepsilon.
\end{equation}

On the other hand, since $x \in X$, we have that $x_je_j\to 0$ in $X$ as $j\to\infty$; with \eqref{eq-cond3} this implies that
\[
e_{q+n_k}= \frac{1}{x_{q+n_k}}x_{q+n_k}e_{q+n_k} \to 0
\]
as $k\to\infty$, and therefore by the continuity of $B$
\begin{equation}\label{eq-app2}
\|e_{p+n}\| = \|e_{q+2p-2q+n_k-n_1}\| = \|B^{n_1-2p+2q}e_{q+n_k}\|<\varepsilon
\end{equation}
whenever $k$ is sufficiently large, where we have noted that $n_1\geq 2p-2q$. 

Altogether we have found some $n\geq 0$ such that \eqref{eq-app1} and \eqref{eq-app2} hold. Thus $B$ is hypercyclic by Theorem \ref{t-hcbil}.
\end{proof}

We have again the following consequence.

\begin{corollary}\label{c-perbil}
Let $X$ be a Fr\'echet sequence space over $\Z$ in which $(e_n)_n$ is a basis. If a weighted backward shift on $X$ has a non-trivial periodic point then it is hypercyclic.
\end{corollary}

The result also allows us to extend a result of Costakis and Parissis \cite[Proposition 5.1]{CoPa12} from $\ell^p(\Z)$ to general Fréchet sequence spaces:

\begin{corollary}\label{c-0-1recbil}
Let $X$ be a Fr\'echet sequence space over $\Z$ in which $(e_n)_n$ is a basis, and let $B_w$ be a weighted backward shift on $X$. Then the following assertions are equivalent:
\begin{enumerate}
\item[\rm (i)] $B_w$ is hypercyclic;
\item[\rm (ii)] $B_w$ is recurrent;
\item[\rm (iii)] $B_w$ has a non-zero recurrent vector.
\end{enumerate}
\end{corollary}

\subsection{Existence of weak sequential limit points}\label{s-wslim} 
We turn to non-zero weak sequential limit points. First of all, for very special spaces, the existence of a non-zero weak sequential limit point does imply hypercyclicity. Recall that a Banach space has the Schur property if every weakly convergent sequence converges; and a Fréchet space is a Montel space if every bounded set is relatively compact. In these spaces, any weakly convergent sequence has (at least) a convergent subsequence. Thus the following is a consequence of Theorem \ref{t-0-1-bilat}.

\begin{theorem}\label{t-SchurMontel}
Let $X$ be a Banach sequence space over $\Z$ with the Schur property or a Fréchet-Montel sequence space over $\Z$. Suppose that $(e_n)_n$ is a basis in $X$, and let $B_w$ be a weighted backward shift on $X$. 

If there exists an orbit with a non-zero weak sequential limit point, then $B_w$ is hypercyclic.
\end{theorem}

This holds, in particular, for the Schur space $\ell^1(\Z)$ and the Fréchet space $H(\C^*)$, see Example \ref{ex-spaces}, which is Montel by the classical Montel theorem. We will show in Corollary \ref{c-weakseqlp} that this result holds also on $\ell^p(\Z)$, for $1\le p<\infty.$

On the contrary, in $c_0(\Z)$, even a large supply of vectors with non-zero weak sequential limit points does not guarantee that a bilateral shift is hypercyclic.  

\begin{example}\label{ex-c0}
Let $X=c_0(\mathbb Z)$ and $B_w$ the weighted backward shift with weights 
\[
w_n= \begin{cases}
     2, &\text{ if } n> 0,\\ 
		 1, &\text{ if } n\leq 0. 
    \end{cases}
\]
Then, $B_w$ is an invertible operator that has a norm dense set of weakly sequentially recurrent vectors but is not hypercyclic. 

Moreover, $B_w$ is weakly hypercyclic but not weakly sequentially hypercyclic.
\end{example}

\begin{proof}
The non-hypercyclicity follows from a well-known characterization, see \cite[Theorem 4.13]{GrPe11}. As for the additional assertions, 
Chan and Sanders \cite[proof of Corollary 3.3]{ChSa04} have shown that $B_w$ is weakly hypercyclic on $\ell^2(\Z)$; since this space is continuously and densely included in $c_0(\Z)$ in the norm sense and therefore also in the weak sense, $B_w$ is also weakly hypercyclic on $c_0(\Z)$. Moreover, by Bès, Chan and Sanders \cite{BCS05}, every weakly sequentially hypercyclic backward shift on $c_0(\Z)$ is hypercyclic.

It remains to show that $B_w$ has a norm dense set of weakly sequentially recurrent vectors. Note that this is much stronger than saying that $B_w$ is weakly sequentially recurrent. We start by showing that there are weakly sequentially recurrent vectors arbitrarily close to $e_0$. To see this, let $(n_k)_{k\geq 1}$ be a sequence of integers with $n_1\geq 1$ and such that, for every $k\geq 2$,
\begin{equation}\label{nj}
n_k\geq\sum_{j=1}^{k-1}n_j+k.
\end{equation} 
We then let
\[
A= \{ n_{j_1}+n_{j_2}+\ldots+n_{j_l}: 1\leq j_1<j_2<\ldots<j_l,\ l\geq 1\}
\]
and 
\[
x=e_0+\sum_{n\in A} \frac{1}{2^n}e_n.
\]
Clearly, $x$ belongs to $c_0(\Z)$, and since $n_1$ is the minimal element of $A$ we have that
\begin{equation}\label{dense}
\|x-e_0\| \leq \frac{2}{2^{n_1}}.
\end{equation}

Since, for any $n,m\geq 0$,
\[
B_w^m\Big(\frac{1}{2^n}e_n\Big) =
\begin{cases}
     e_{n-m}, &\text{ if } m\geq n,\\ 
		 \frac{1}{2^{n-m}}e_{n-m}, &\text{ if } m<n, 
    \end{cases}
\]
we find that, for any $k\geq 1$,
\begin{equation}\label{bwnk}
B_w^{n_k} x = e_{-n_k}+\sum_{\substack{n\in A\\n<n_k}} e_{n-n_k} + e_0 +\sum_{\substack{n\in A\\n_k<n<n_{k+1}}} \frac{1}{2^{n-n_k}}e_{n-n_k}+\sum_{\substack{n\in A\\n\geq n_{k+1}}} \frac{1}{2^{n-n_k}}e_{n-n_k};
\end{equation}
note that $n_k\in A$.

First, if $n\in A$, $n<n_k$, then we have in view of \eqref{nj} that there are $j_1<j_2<\ldots<j_l<k$ with
\[
n=n_{j_1}+n_{j_2}+\ldots+n_{j_l}\leq \sum_{j=1}^{k-1}n_j\leq n_k-k,
\]
so that $n-n_k\leq -k$. Thus the first two terms of the right-hand side in \eqref{bwnk} tend weakly to 0, because, for any $y\in \ell^1(\mathbb Z)$,
\[
\Big|\Big\langle e_{-n_k}+\sum_{\substack{n\in A\\n<n_k}} e_{n-n_k}, y\Big\rangle\Big|\leq \sum_{n\leq -k}|y_n|\to 0.
\]
Next, \eqref{nj} also implies that $n\in A$, $n_k<n<n_{k+1}$ if and only if there are $j_1<j_2<\ldots<j_l<k$ such that
\[
n=n_{j_1}+n_{j_2}+\ldots+n_{j_l}+n_k,
\]
so that the corresponding numbers $m=n-n_k$ cover exactly the elements in $A$ that are less than $n_k$. Consequently,
\[
e_0 +\sum_{\substack{n\in A\\n_k<n<n_{k+1}}} \frac{1}{2^{n-n_k}}e_{n-n_k} = e_0 + \sum_{\substack{m\in A\\m<n_k}} \frac{1}{2^{m}}e_{m}\to x
\]
in norm and hence also weakly.

Finally, if $n\in A$, $n\geq n_{k+1}$, then by \eqref{nj}
\[
n-n_k\geq n_{k+1}-n_k\geq k,
\]
which implies that
\[
\sum_{\substack{n\in A\\n\geq n_{k+1}}} \frac{1}{2^{n-n_k}}e_{n-n_k}\to 0,
\]
again in norm and hence weakly.

Altogether we have found that
\begin{equation}\label{bwweak}
B_w^{n_k}x\to x
\end{equation}
in the weak topology.

To finish the proof we use the fact that $B_w$ is invertible. Given any scalars $a_n$, $n=-N,\ldots,N$, $N\geq 1$, the (weak) continuity of $B_w$ implies with \eqref{bwweak} that
\[
B_w^{n_k}\Big(\sum_{n=-N}^N a_n B_w^n x\Big)=\sum_{n=-N}^N a_n B_w^n (B_w^{n_k} x) \to \sum_{n=-N}^N a_n B_w^n x
\]
in the weak topology, hence $\sum_{n=-N}^N a_n B_w^n x$ is weakly sequentially recurrent. Moreover, we have by \eqref{dense} that
\[
\Big\|\sum_{n=-N}^N a_n B_w^n x - \sum_{n=-N}^N a_n B_w^n e_0 \Big\|\leq \frac{2}{2^{n_1}}\Big\|\sum_{n=-N}^N a_n B_w^n\Big\|.
\]
Now note that any finite sequence can be written in the form $\sum_{n=-N}^N a_n B_w^n e_0$ for some $N\geq 1$ and some scalars $a_n$. Since one may choose $n_1$ arbitrarily large, the density of the finite sequences in $c_0(\Z)$ then implies that the weakly sequentially recurrent vectors for $B_w$ form a dense set.
\end{proof}

A modification of the previous example shows that the same phenomenon can occur in a reflexive Banach sequence space.

\begin{example}\label{ex-refl}
There is a reflexive Banach sequence space $X$ in which $(e_n)_{n\in\Z}$ is an unconditional basis and an invertible weighted backward shift $B_w$ on $X$ that has a norm dense set of weakly sequentially recurrent vectors but is not hypercyclic. 

Moreover, $B_w$ is weakly hypercyclic but not weakly sequentially hypercyclic.
\end{example}

\begin{proof}
We consider the same weights as in the previous example,
\[
w_n= \begin{cases}
     2, &\text{ if } n> 0,\\ 
		 1, &\text{ if } n\leq 0. 
    \end{cases}
\]
We demand here slightly more of the sequence $(n_k)_{k\geq 1}$: let $n_1\geq 2$ and, for every $k\geq 2$,
\[
n_k\geq\sum_{j=1}^{k-1}n_j+n_{k-1}+k.
\]
Let again
\[
A= \{ n_{j_1}+n_{j_2}+\ldots+n_{j_l}: 1\leq j_1<j_2<\ldots<j_l,\ l\geq 1\}.
\]
It then follows that, for any $n\in A$ with $n<n_k$, $k\geq 1$, 
\begin{equation}\label{nj'}
-n_k < n-n_k < -n_{k-1}.
\end{equation} 
We consider the sets of negative integers
\[
I_k = [-n_k,-n_{k-1}),\ \ k\geq 1
\]
with $n_0=0$, and the numbers $p_k:=\text{card}(I_k)=n_k-n_{k-1}\geq 2$, $k\geq 1$.

The idea of our construction is to replace the space $c_0(\Z)$ of Example \ref{ex-c0} by a smaller reflexive Banach sequence space $X$ which is big enough so that the argument of the example continues to hold. With this aim, we regard the space $X$ of all sequences $x=(x_n)_{n\in\Z}$ such that
\[
\|(x_n)_{n\geq 0}\|^2_2 + \sum_{k=1}^\infty \|(x_n)_{n\in I_k}\|_{p_k}^2 <\infty,
\]
with the canonical norm. Then $X$ is reflexive as the direct $\ell^2$-sum of the reflexive spaces $\ell^2(\N_0)$ and $\ell^{p_k}(I_k)$, $k\geq 1$, and $(e_n)_n$ is an unconditional basis. 

By considering separately the sequences $x\in X$ with support inside and outside $\{-n_k:k\geq 0\}$ one sees that the weighted backward shift $B_w$ maps $X$ into $X$ and is therefore continuous; similarly, $B_w$ is seen to be invertible. And since $\|e_n\|=1$ for all $n\in \Z$, $B_w$ cannot be hypercyclic in view of \cite[Theorem 4.13]{GrPe11}.

To see that it is not even weakly sequentially hypercyclic, note that $X$ is continuously and densely included in $c_0(\Z)$ in the norm sense and therefore also in the weak sense. Thus, if $B_w$ were weakly sequentially hypercyclic on $X$, it would also be so on $c_0(\Z)$, which is not the case by the previous example. By the same argument, since $B_w$ is weakly hypercyclic on $\ell^2(\Z)$ by \cite[Corollary 10.28]{BaMa09}, it is weakly hypercyclic on the larger space $X$ that contains $\ell^2(\Z)$ as a dense subspace.

Now, let again
\[
x=e_0+\sum_{n\in A} \frac{1}{2^n}e_n.
\]
Then $x$ belongs to $X$, and we have that
\begin{equation}\label{dense'}
\|x-e_0\| \leq \frac{1}{2^{n_1}}\frac{2}{\sqrt{3}}.
\end{equation}
The identity \eqref{bwnk} for $B_w^{n_k}x$, $k\geq 1$, continues to hold, and the sum of the last three terms on the right-hand side converges to $x$ in norm, and hence weakly: replacing the norm of $c_0(\N)$ by that of $\ell^2(\N)$ preserves the validity of our previous argument.

Let us now show that 
\[
e_{-n_k}+\sum_{\substack{n\in A\\n<n_k}} e_{n-n_k}\to 0
\]
weakly in $X$ as $k\to\infty$. To see this, let $y=(y_n)_{n\in\Z}\in X^*$. Then 
\[
\|(y_n)_{n\geq 0}\|^2_2 + \sum_{k=1}^\infty \|(y_n)_{n\in I_k}\|_{p^*_k}^2 <\infty,
\]
where $\frac{1}{p_k^*}+ \frac{1}{p_k}=1$, $k\geq 1$. We have by \eqref{nj'} that for any $n\in A$, $n<n_k$, $k\geq 1$,
\[
n-n_k\in I_k,
\]
and also $-n_k\in I_k$. Thus
\[
\Big|\Big\langle e_{-n_k}+\sum_{\substack{n\in A\\n<n_k}} e_{n-n_k}, y\Big\rangle\Big|\leq \sum_{n\in I_k}|y_n|\leq \text{card}(I_k)^{1/p_k}\|(y_n)_{n\in I_k}\|_{p^*_k}\to 0;
\]
note that $p_k=\text{card}(I_k)$. This establishes the previous claim.

Altogether we have again that
\[
B_w^{n_k}x\to x
\]
weakly in $X$. 

From there we conclude the proof with \eqref{dense'} as for Example \ref{ex-c0}.
\end{proof}

While the existence of a norm dense set of weakly sequentially recurrent vectors need not imply hypercyclicity for bilateral shifts, one may ask whether the existence of a non-zero weakly sequentially recurrent vector implies that the set of such vectors is necessarily norm dense. We will give a positive answer for many sequence spaces, including $c_0(\Z)$. 

To motivate what follows, note that one of the main advantages of working with unconditional bases is that if $(x^{(k)})_k$ is a sequence of vectors that converges to 0 in the space $X$ and if another sequence $(y^{(k)})_k$ satisfies that, for some constant $K>0$,
\[
|y^{(k)}_n|\leq K |x^{(k)}_n|,\ n\in\Z
\]
then $(y^{(k)})_k$ also converges to 0. We will say that $(y^{(k)})_k$ is \textit{coordinatewise controlled} by $(x^{(k)})_k$. Now, what we will need is a similar property for weak sequential convergence.

\begin{definition}
We say that that a Fréchet sequence space $X$ has property (BW) if, for any two sequences $(x^{(k)})_k$ and $(y^{(k)})_k$ of vectors in $X$, if $(x^{(k)})_k$ converges weakly to $0$ and $(y^{(k)})_k$ is coordinatewise controlled by $(x^{(k)})_k$, then $(y^{(k)})_k$ converges weakly to $0$.
\end{definition}

Many of the standard Banach sequence spaces have property (BW). Recall that a basis $(e_n)_n$ in a Banach sequence space is shrinking provided that the sequence of coordinate functionals $(e_n^*)_n$ is a basis for $X^*$. In particular, any basis in a reflexive Banach space is shrinking, see \cite[Chapter I, Corollary 12.2]{Sin70}.

\begin{proposition}\label{p-bw}
Let $X$ be a Banach sequence space over $\Z$ in which $(e_n)_n$ is an unconditional basis. If $(e_n)_n$ is shrinking, in particular if $X$ is a reflexive Banach space, then $X$ satisfies property \emph{(BW)}.
\end{proposition}

\begin{proof}
Let $(x^{(k)})_k$ and $(y^{(k)})_k$ be sequences of vectors in $X$ such that $(x^{(k)})_k$ converges weakly to 0 and that $(y^{(k)})_k$ is coordinatewise controlled by $(x^{(k)})_k$. Note that then $(y^{(k)})_k$ is coordinatewise convergent to 0.

Since $(x^{(k)})_k$ is weakly bounded, hence bounded, then so is $(y^{(k)})_k$ by unconditionality of the basis. Now the fact that $(e_n)_n$ is shrinking implies that $(y^{(k)})_k$ converges weakly to $0$, see \cite[Chapter II, Theorem 4.2]{Sin70}.
\end{proof}

Thus, in particular, $c_0(\Z)$ and $\ell^p(\Z)$, $1< p<\infty$, have property (BW). Note that $\ell^1(\Z)$ also has property (BW) by the Schur property, but that space has already been dealt with above, see Theorem \ref{t-SchurMontel}.

We will also need the following.

\begin{lemma}\label{l-bw}
Let $X$ be a Fréchet sequence space over $\Z$ in which $(e_n)_n$ is a basis with property \emph{(BW)}. Let $(x^{(k)})_k$ and $(y^{(k)})_k$ be sequences of vectors in $X$ so that  $(x^{(k)})_k$ converges weakly, $(y^{(k)})_k$ is coordinatewise controlled by $(x^{(k)})_k$ and, for every $n\in \Z$, $(y^{(k)}_n)_k$ is eventually $0$. Then $(y^{(k)})_k$ converges weakly to $0$.
\end{lemma}

\begin{proof} Let $x\in X$ be the weak limit of $(x^{(k)})_k$. Since we also have coordinatewise convergence, there is a sequence $(N_k)_k$ of integers with $N_k\to\infty$ such that $\sum_{|n|\leq N_k} \|(x_n-x_n^{(k)})e_n\|\to 0$. By decreasing the $N_k$, if necessary, we can assume by hypothesis that $y^{(k)}$ has support in  $[-N_k,N_k]^c$, $k\geq 1$.

We consider $z^{(k)}:=\sum_{ |n|> N_k} x_n^{(k)}e_n$. Since
\[
z^{(k)} = (x^{(k)}-x) + \sum_{ |n|> N_k} x_ne_n + \sum_{|n|\leq N_k} (x_n-x_n^{(k)})e_n,
\]
the basis assumption and the definition of $(N_k)_k$ imply that $(z^{(k)})_k$ converges weakly to 0 and that it dominates $(y^{(k)})_k$ coordinatewise. By property (BW), $(y^{(k)})_k$ converges weakly to 0.
\end{proof}

We can now deduce our main result on weak sequential recurrence.

\begin{theorem}\label{t-weaksr}
Let $X$ be a Fr\'echet sequence space over $\Z$ in which $(e_n)_n$ is a basis with property \emph{(BW)}, and let $B_w$ be a weighted backward shift on $X$. 

If $B_w$ has a non-zero weakly sequentially recurrent vector, then it has a norm dense set of such vectors.
\end{theorem}

\begin{proof}
By conjugacy, which also preserves property (BW), it suffices to consider the unweighted shift $B$.

Let $y\in X$ be of the form $y=\sum_{j=-p}^p y_ne_n$ with $p\geq 0$, and let $\varepsilon >0$. The continuity of the operator
\[ 
T:=\sum_{n=0}^{2p} y_{p-n} B^n
\]
with $Te_p=y$ implies that there is some $\eta>0$ such that, for all $z\in X$,
\begin{equation}\label{eq-T}
\|z-e_p\|<\eta \ \Longrightarrow \ \|Tz-y\|<\varepsilon.
\end{equation}

For the sequel, we fix these $p$ and $\eta$. 

Now, by assumption, there is a weakly sequentially recurrent vector $x\neq 0$ for $B$. By applying $B$ sufficiently often and noting that $B$ is weakly continuous, we can assume that 
\begin{equation}\label{eq-xr}
\delta:=|x_p|>0.
\end{equation}

There is then an infinite set $D\subset\N_0$ such that
\[
(B^nx)_{n\in D}
\]
converges weakly to $x$.

Let $q\geq 1$. By weak recurrence, which implies coordinatewise recurrence, there is some $m_q\in D$, $m_q\geq q$, such that
\[
|x_{p+m_q}-x_p|=|[B^{m_q}x-x]_p|<\frac{\delta}{2},
\]
and hence $|x_{p+m_q}|>\frac{\delta}{2}$.

Again by coordinatewise recurrence, there is an infinite set $D_q\subset D$ such that, for any $n\in D_q$, 
\begin{equation}\label{eq-abs}
|[B^nx-x]_{p+j}|<\frac{\delta}{2^{q+1}},\ \ 0\leq j\leq m_q.
\end{equation}
In particular we obtain that 
\[
|x_{p+n+m_q}-x_{p+m_q}|= |[B^nx-x]_{p+m_q}|< \frac{\delta}{4}
\]
and thus $|x_{p+n+m_q}|>\frac{\delta}{4}$. Since $x=\sum_{j\in\Z} x_je_j$ converges in $X$, we have that $x_je_j\to 0$ as $j\to\infty$ and hence
\[
e_{p+n+m_q}= \frac{1}{x_{p+n+m_q}}x_{p+n+m_q}e_{p+n+m_q}\to 0, \ n\to\infty, n\in D_q.
\]
Applying $B$ up to $2m_q$ times and noting that $q\leq m_q$, we find that for all $n\in D_q$ sufficiently large,
\begin{equation}\label{eq-norme}
\|e_{p+n+j}\| < \frac{1}{q2^{q+1}}\eta, \ \ |j|\leq q.
\end{equation}
We may assume without loss of generality that this is the case for all $n\in D_q$. 

We now define inductively a sequence of integers $(n_k)_{k\geq 0}$ by taking $n_0=0$ and then, for every $k\geq 1$,
\begin{align}
N_{k-1}&:= \sum_{j=0}^{k-1}n_j,\label{nj1}\\
q_k&:= 2N_{k-1}+1,\label{nj2}\\
n_{k}&\in D_{q_k}\ \text{and}\ n_{k}\geq N_{k-1}+k.\label{nj4}
\end{align} 
We let
\[
A= \{ n_{j_1}+n_{j_2}+\ldots+n_{j_m}: 1\leq j_1<j_2<\ldots<j_l,\ l\geq 1\}
\]
and define
\begin{equation}\label{eq-x}
z=e_p+\sum_{n\in A} e_{p+n}.
\end{equation}

Our first task is to show that this series converges. We have that
\[
\sum_{n\in A} \|e_{p+n}\| = \sum_{k=1}^\infty \sum_{\substack{n\in A\\n_k\leq n<n_{k+1}}} \|e_{p+n}\|.
\]
Now, for $k\geq 1$ and $n\in A$ with $n_k\leq n<n_{k+1}$, we have by \eqref{nj1} and \eqref{nj4} that $0\leq j:=n-n_k\leq N_{k-1}$, and hence by \eqref{eq-norme}, \eqref{nj2} and \eqref{nj4} that
\begin{align*}
\sum_{\substack{n\in A\\n_k\leq n<n_{k+1}}}\| e_{p+n}\|&=\sum_{\substack{n\in A\\n_k\leq n<n_{k+1}}}\| e_{p+n_k+(n-n_k)}\| \leq \sum_{0\leq j\leq N_{k-1}}\frac{1}{q_k2^{q_k+1}}\eta\\
&\leq \sum_{0\leq j\leq N_{k-1}}\frac{1}{(2N_{k-1}+1)2^{k}}\eta = \frac{1}{2^{k}}\eta,
\end{align*}
which implies the convergence of \eqref{eq-x}, so that $z$ is a well-defined vector in $X$. Moreover,
\begin{equation}\label{eq-T2}
\|z-e_p\|\leq \eta.
\end{equation}

Next, for any $k\geq 1$, we have that
\begin{equation}\label{bwnk2}
B^{n_k} z = e_{p-n_k}+\sum_{\substack{n\in A\\n<n_k}} e_{p+n-n_k} + e_p +\sum_{\substack{n\in A\\n_k<n<n_{k+1}}} e_{p+n-n_k}+\sum_{j=k+1}^\infty \sum_{\substack{n\in A\\n_j\leq n < n_{j+1}}} e_{p+n-n_k}.
\end{equation}

Let's look at the right-hand side of this identity.

First, we have seen in the proof of Example \ref{ex-c0} that, for any $k\geq 1$, $\{n-n_k: n \in A, n_k<n<n_{k+1}\} = \{n \in A, n<n_{k}\}$, so that we have
\begin{equation}\label{eq-conv1}
e_p +\sum_{\substack{n\in A\\n_k<n<n_{k+1}}}e_{p+n-n_k} = e_p +\sum_{\substack{n\in A\\n<n_{k}}} e_{p+n}\to z
\end{equation}
in norm and hence weakly.

Secondly, for any $j\geq k+1$,
\[
\sum_{\substack{n\in A\\n_j\leq n < n_{j+1}}} \|e_{p+n-n_k}\| = \sum_{\substack{n\in A\\n_j\leq n < n_{j+1}}} \|e_{p+n_j + (n-n_j -n_k)}\|.
\]
Since, for any $n\in A$ with $n_j\leq n < n_{j+1}$ we have by \eqref{nj1} that
\[
-N_{j-1}\leq -n_k\leq m:=n-n_j-n_k\leq n-n_j\leq N_{j-1},
\]
we obtain with \eqref{eq-norme}, \eqref{nj2} and \eqref{nj4} that
\[
\sum_{\substack{n\in A\\n_j\leq n < n_{j+1}}} \|e_{p+n-n_k}\| \leq  \sum_{\substack{n\in A\\n_j\leq n < n_{j+1}}} \frac{1}{q_{j}2^{q_{j}+1}}\eta\leq \sum_{|m|\leq N_{j-1}}\frac{1}{(2N_{j-1}+1)2^{j}}\eta = \frac{1}{2^{j}}\eta.
\]
Thus
\begin{equation}\label{eq-conv2}
\Big\|\sum_{j=k+1}^\infty \sum_{\substack{n\in A\\n_j\leq n < n_{j+1}}} e_{p+n-n_k}\Big\|\leq \frac{1}{2^{k}}\eta \to 0.
\end{equation}

Thirdly, let $x$ be the initial weakly sequentially recurrent vector. If $n\in A$, $n<n_{k}$, then there are $l\geq 1$ and $1\leq j_1<j_2<\ldots<j_l\leq k-1$ with
\[
n=n_{j_1}+n_{j_2}+\ldots+n_{j_{l-1}}+n_{j_l}.
\]
Let $m=1,\ldots,l$. Since $n_{j_m} \in D_{q_{j_m}}$ and, for $m\geq 2$,
\[
0\leq n_{j_1}+\ldots + n_{j_{m-1}} \leq N_{j_{m-1}}\leq N_{j_{m}-1}\leq q_{j_{m}}\leq m_{q_{j_m}},
\]
we have by \eqref{eq-abs} that, for any $m\geq 1$,
\begin{align*}
|x_{p+ n_{j_1}+\ldots+n_{j_{m-1}}+n_{j_m}} - x_{p+ n_{j_1}+\ldots+n_{j_{m-1}}}| &= |[B^{n_{j_m}}x-x]_{p+n_{j_1}+\ldots+n_{j_{m-1}}}|\\
&\leq \frac{\delta}{2^{q_{j_m}+1}}\leq \frac{\delta}{2^{m+1}}, 
\end{align*}
where an empty sum is 0. Telescoping we find that
\[
|x_{p+n}-x_p|\leq \sum_{m=1}^l \frac{\delta}{2^{m+1}} \leq \frac{\delta}{2}
\]
and therefore by \eqref{eq-xr}
\begin{equation}\label{eq-xrn}
|x_{p+n}| \geq \frac{\delta}{2}.
\end{equation}

When we now define
\[
y^{(k)} := e_p + \sum_{\substack{n\in A\\n<n_{k}}} e_{p+n}, \ \ k\geq 1,
\]
then it follows from \eqref{eq-xrn} that the sequence $(y^{(k)})_{k\geq 1}$ is coordinatewise controlled by $(x)_k$, a constant sequence of vectors in $X$. Since the operator $B$ is a simple shift, we deduce that the sequence $(B^{n_k}y^{(k)})_{k\geq 1}$ is coordinatewise controlled by the sequence $(B^{n_k}x)_{k\geq 1}$. Now, since each $n_k$ belongs to $D$, $(B^{n_k}x)_{k\geq 1}$ is weakly convergent. And finally, since
\[
B^{n_k}y^{(k)}=e_{p-n_k} + \sum_{\substack{n\in A\\n<n_{k}}} e_{p+n-n_k}
\]
and, for any $n\in A$, $n<n_{k}$, by \eqref{nj4},
\[
p+n-n_k \leq p+N_{k-1}-n_k \leq p-k,
\]
we see that, for any $j\in\Z$, $([B^{n_k}y^{(k)}]_j)_{k\geq 1}$ is eventually 0. By property (BW), via Lemma \ref{l-bw}, this implies that 
\begin{equation}\label{eq-conv3}
e_{p-n_k} + \sum_{\substack{n\in A\\n<n_{k}}} e_{p+n-n_k}\to 0
\end{equation}
in the weak topology.

Altogether, we see from \eqref{bwnk2}, \eqref{eq-conv1}, \eqref{eq-conv2} and \eqref{eq-conv3} that
\[
B^{n_k}z\to z
\]
in the weak topology.

To finish the proof, we apply the (weakly) continuous operator $T=\sum_{n=0}^{2p} y_{p-n} B^n$. Since it commutes with $B$, we have that
\[
B^{n_k}Tz\to Tz
\]
in the weak topology, and by \eqref{eq-T} and \eqref{eq-T2} that
\[
\|Tz-y\|<\varepsilon.
\]
We have therefore shown that in every neighbourhood of any finite sequence there exists a weakly sequentially recurrent vector. Since the finite sequences are dense in $X$, the theorem is proved.
\end{proof}

Apart from spaces with very strong properties, see Theorem \ref{t-SchurMontel}, we have not been able to obtain a zero-one law for weak sequential limit points. In particular, the following might be of interest.

\begin{question}\label{q-weakseqrec}
Let $B_w$ be a weighted backward shift on a reflexive Banach sequence space over $\Z$ in which $(e_n)_n$ is a basis. Does the existence of a non-zero weak sequential limit point imply weak sequential recurrence or even the existence of a norm dense set of weakly sequentially recurrent vectors?
\end{question}

Note that by Example \ref{ex-refl} we cannot expect $B_w$ to be hypercyclic. Below, we will get a positive answer for the spaces $\ell^p(\Z)$, $1\leq p<\infty$, see Corollary \ref{c-weakseqlp}, in which case $B_w$ even turns out to be hypercyclic.

\subsection{Existence of weak*-sequential limit points}\label{s-w*slim}
Recall that if $X$ is a Banach sequence space in which $(e_n)_n$ is a boundedly complete basis then it is isomorphic to the dual of the Banach sequence space $Z$ in which $(e_n^*)_n$ is a basis, see \cite[Theorem 3.2.10]{AlKa06}.
In this case we may look at weak*-sequential limit points. We have the following analogue of Theorem \ref{t-weaksr}.

\begin{theorem}\label{t-weak*sr}
Let $X$ be a Banach sequence space over $\Z$ in which $(e_n)_n$ is an unconditional boundedly complete basis, and let $B_w$ be a weighted backward shift operator on $X$. 

If $B_w$ has a non-zero weak*-sequentially recurrent vector, then it has a norm dense set of such vectors.
\end{theorem}

\begin{proof} 
Note first that since $B_w^*(e_n^*)=w_{n+1}e^*_{n+1}\in Z$ it follows that  $B_w^*|_Z:Z\to Z$ is a well defined operator. This implies that $B_w=(B_w^*|_{Z})^*$ and hence $B_w$ is weak*-continuous. Also, weak*-convergence implies coordinatewise convergence.  

Let us now show that the analogue of Lemma \ref{l-bw} holds. Thus let $(x^{(k)})_k$ and $(y^{(k)})_k$ be sequences of vectors in $X$ so that $(x^{(k)})_k$ converges weak*, $(y^{(k)})_k$ is coordinatewise controlled by $(x^{(k)})_k$ and, for every $n\in \Z$, $(y^{(k)}_n)_k$ is eventually 0. It follows from the Banach-Steinhaus theorem that $(x^{(k)})_k$ is bounded in $X$, and hence also $(y^{(k)})_k$ by unconditionality of the basis in $X$.

Now let $z\in Z$ and $\varepsilon>0$. Since $(e_n^*)_n$ is a basis of $Z$, there is some $N\in\N$ such that
\[
\Big\|\sum_{|n|\geq N} z_ne_n^*\Big\|<\frac{\varepsilon}{1+ \sup_{k\geq 1}\| y^{(k)}\|}.
\]
Since $y^{(k)}\to 0$ coordinatewise, there is then some $K\geq 1$ such that, for any $k\geq K$,
\[
|\langle y^{(k)},z\rangle|\leq \Big|\sum_{|n|< N} y^{(k)}_nz_n\Big| + \| y^{(k)}\| \Big\|\sum_{|n|\geq N} z_ne_n^*\Big\|< 2\varepsilon.
\]
Thus, $(y^{(k)})_k$ converges weak* to $0$.

Now the proof can be done exactly as that of Theorem \ref{t-weaksr}: it suffices to replace throughout ``weak'' and ``weakly'' by ``weak*-''. 
\end{proof}

\subsection{Existence of weak and weak* limit points}\label{s-weaklim} 
As we have mentioned above, Shkarin \cite{Shk07} has shown that every weakly hypercyclic weighted backward shift on $\ell^p(\Z)$, $1\leq p< 2$, is hypercyclic, see also \cite[Proposition 10.20]{BaMa09}. By a slight modification of the proof we obtain the following zero-one law.

\begin{theorem}\label{t-0-1 bil ell_p weak p<2}
Let $B_w$ be a weighted backward shift on $\ell^p(\mathbb Z)$, $1\leq p<2$. If $B_w$ admits an orbit with a weak limit point $y$ such that $\sup\{n\geq 0: y_n\neq 0\}=\infty$ then $B_w$ is hypercyclic.
\end{theorem}

\begin{proof}
Suppose that $B_w$ is not hypercyclic. Then there are some $q\geq 0$ and $\varepsilon>0$ such that, for all $n\geq 0$, 
\begin{equation}\label{eq: 0-1 weak p<2}
\max\big(|w_{q+1} \cdots w_{q+n}|^{-1} , |w_{q-n+1}\cdots w_{q}|\big)\geq \varepsilon 
\end{equation}
by Theorem \ref{t-hcbil} and using the conjugacy explained in the Introduction.

Assume that $y$ with $\sup\{n\geq 0: y_n\neq 0\}=\infty$ is a weak limit point of the orbit of some vector $x$. By weak continuity of $B_w$
we may assume that $y_q\neq 0$, and then that $y_q=2$.

Moreover, suppose that there is some $n\geq 0$ such that $B_w^nx=y$. If we also had $B_w^mx=y$ for some $m>n$, $B_w^nx$ would be a non-zero periodic point, which implies that $B_w$ is hypercyclic by Corollary \ref{c-perbil}. Thus, replacing $x$ by $B_w^{n+1}x$, we see that we can assume that $y$ does not belong to the orbit of $x$.

Define $\Lambda =\{n\geq 0:  |[B_w^{n}x]_q|>1\}$.  Since $y$ is a weak limit point of $(B_w^{n}x)_{n}$, $\Lambda$ is an infinite set. Moreover, $\{B_w^{n}x:n\in\Lambda\}$ is not weakly closed because its weak closure contains $y$, which is outside this set.

Now we can follow the proof of \cite[Proposition 10.20]{BaMa09} to show that $\{B_w^{n}x:n\in\Lambda\}$ is weakly closed, thus arriving at a contradiction.

Indeed, for any $n\in\Lambda$ we have that 
\begin{align*}
\|B_w^{n}x\|^p & \ge \sum_{m\in\Lambda} |w_{q+m-n+1} \cdots w_{q+m}|^{p}|x_{q+m}|^p    \\
& = \sum_{m\in\Lambda, m\le n} |w_{q+m-n+1} \cdots w_{q}|^{p}|[B_w^{m}x]_q|^p + \sum_{m\in\Lambda, m> n} |w_{q+1} \cdots w_{q+m-n}|^{-p}|[B_w^{m}x]_q|^p \\
& \ge \sum_{m\in\Lambda, m\le n} |w_{q+m-n+1} \cdots w_{q}|^{p} + \sum_{m\in\Lambda, m> n} |w_{q+1} \cdots w_{q+m-n}|^{-p}.
\end{align*}
Define the infinite matrix $(c_{n,m})_{n,m\in \Lambda}$ by $c_{n,m}=|w_{q+m-n+1} \cdots w_{q}|^{p}$ if $ m\le n$ and $c_{n,m}=|w_{q+1} \cdots w_{q+m-n}|^{-p}$ if $m> n$. Then by \eqref{eq: 0-1 weak p<2} we can apply \cite[Lemma 10.22]{BaMa09} to conclude that $\sum_{n\in\Lambda}(\sum_{m\in\Lambda}c_{n,m})^{-r}<\infty$ for any $r>1$, and thus $\sum_{n\in\Lambda}\|B_w^{n}x\|^{-rp}<\infty$. Taking $r>1$ such that $rp<2$,  \cite[Corollary 10.7]{BaMa09} implies that $\{B_w^{n}x:n\in \Lambda\}$ is weakly closed, which is what we need to show.
\end{proof}

Since no non-zero weakly recurrent vector can be supported in $(-\infty,N]$ for some $N\in\Z$, we obtain the following.

\begin{corollary}\label{c-0-1 bil ell_p}
Let $B_w$ be a weighted backward shift on $\ell^p(\mathbb Z)$, $1\leq p<2$. If $B_w$ has a non-zero weakly recurrent vector, then it is hypercyclic.
\end{corollary}

As we have already mentioned, the results fail on any space $\ell^p(\Z)$, $p\geq 2$, because, by Chan and Sanders \cite{ChSa04}, they support a weakly hypercyclic backward shift that is not hypercyclic. The same is then true on $c_0(\Z)$ because weak hypercyclicity passes from $\ell^2(\Z)$ to the larger space $c_0(\Z)$, and both spaces have the same continuous shifts and the same hypercyclic shifts.

Now, since the spaces $\ell^p(\Z)$, $1<p<2$, are reflexive, the weak*-versions of Theorem \ref{t-0-1 bil ell_p weak p<2} and its corollary hold as well, which is less clear for $\ell^1(\Z)$. In particular, the following seems to be open.

\begin{question}
Is every weak*-hypercyclic weighted backward shift on $\ell^1(\Z)$ hypercyclic?
\end{question}

\subsection{Existence of coordinatewise limit points}\label{s-coordlim} 
The proof of Chan and Seceleanu \cite{ChSe12} of the case $X=\ell^p(\Z)$, $1\leq p<\infty$, of Theorem \ref{t-0-1-bilat} is very different from ours. In particular, they use at a crucial point that the suborbit under consideration is bounded, which we haven't used explicitly above. An analysis of their proof gives in fact a zero-one law for coordinatewise limit points on a large class of spaces.

It will be convenient in this subsection to work with an increasing sequence $(\|\cdot\|_r)_r$ of seminorms that defines the topology of a Fréchet sequence space.

\begin{definition} Let $(X,(\|\cdot\|_r)_r)$ be a Fréchet sequence space that contains all sequences of finite support. Then $X$ is said to have property (BC) if, for any $r\geq 1$ and $\varepsilon >0$, there are $M\in\Z$ and $N\in\N$ such that, for any set $I$ of integers with $I\subset (-\infty, M]$ and $\text{card}(I)=N$, and for any scalars $x_n$, $n\in I$,
\[
\Big\|\sum_{n\in I} x_ne_n\Big\|_r\leq 1 \ \Longrightarrow \ \exists n\in I, \|x_ne_n\|_r < \varepsilon.
\]
\end{definition}

\begin{example}\label{ex-bc}
The spaces $\ell^p(\Z)$, $1\leq p<\infty$, and the space $H(\C^*)$ of Example \ref{ex-spaces} have property (BC), while $c_0(\Z)$ does not. 
\end{example}

\begin{theorem}\label{t-coobil01}
Let $X$ be a Fréchet sequence space over $\Z$ in which $(e_n)_{n\in\Z}$ is an unconditional basis with property \emph{(BC)}, and let $B_w$ be a weighted backward shift on $X$. 

Suppose that there are $x,y\in X$ with $y\neq 0$ and $(n_k)_k$ such that
\[
\{B_w^{n_k}x:k\geq 1\} \text{ is bounded and } B_w^{n_k}x\to y \text{ in } \K^{\Z}.
\]
Then $B_w$ is hypercyclic.
\end{theorem}

\begin{proof}
By conjugacy, which also preserves property (BC), it suffices to consider the unweighted shift $B$.

Suppose that there are $x,y\in X$ with $y\neq 0$ such that the hypothesis holds. Then there is some $q\in\Z$ such that $y_q\neq 0$, where we can assume that $y_q=2$. Also, by the second part of the hypothesis, we can assume that, for all $k\geq 1$,
\begin{equation}\label{eq-xk}
|x_{q+n_k}|=|[B^{n_k}x]_q|\geq 1.
\end{equation}

The remainder of the proof is in two stages. We first show the following claim: for any $r\geq 1$, $\varepsilon>0$ and $m\geq 0$ there are $l\geq 1$ and $k\geq l$ such that
\begin{equation}\label{eq-xk1}
\|e_{q+n_l-n_k}\|_r<\varepsilon \quad \text{and}\quad \|e_{q+n_k-n_l+m}\|_r<\varepsilon.
\end{equation}

To see this, note that by unconditionality of the basis there is some $s\geq 1$ and a constant $C>0$ such that, whenever $z\in X$ and $|a_n|\leq 1$ for all $n\in \Z$, then
\begin{equation}\label{eq-uncond}
\Big\|\sum_{n\in\Z} a_nz_ne_n\Big\|_r\leq C\|z\|_s.
\end{equation}
Now let $K=\sup_{k\geq 1}\|B^{n_k}x\|_s$, and let $N$ and $M$ be given by property (BC), applied to $r$ and $\frac{\varepsilon}{CK}$. We first have, for any $k\geq 1$,
\begin{align*}
1 &\geq \frac{1}{K}\|  B^{n_k} x\|_s = \frac{1}{K}\Big\|  \sum_{n\in\Z} x_n e_{n-n_k}\Big\|_s= \frac{1}{K}\Big\|  \sum_{n\in\Z} x_{q+n} e_{q+n-n_k}\Big\|_s\geq  \Big\| \frac{1}{CK} \sum_{l=m+1}^{m+N}  e_{q+n_l-n_k}\Big\|_r,
\end{align*}
where the last inequality follows from \eqref{eq-xk} and \eqref{eq-uncond}. Then (BC) tells us that for any $k\geq 1$ with $n_k\geq q+n_{m+N}-M$ there is some $l_k\in[m+1, m+N]$ such that
\[
\| e_{q+n_{l_k}-n_k}\|_r<\varepsilon.
\]

On the other hand, since $x\in X$, we have that $x_{q+n_k}e_{q+n_k}\to 0$ as $k\to \infty$, so that $e_{q+n_k}\to 0$ by \eqref{eq-xk} and thus $B^{n_l-m}e_{q+n_k}\to 0$ for any $l\in[m+1,m+N]$; note that $n_l-m\geq 0$. Then there is some $k\geq m+N$ such that
\[
\|e_{q+n_k-n_{l_k}+m}\|_r=\| B^{n_{l_k}-m}e_{q+n_k}\|_r<\varepsilon.
\]
Since $k\geq m+N\geq l_k$, the claim has been shown.

To finish the proof, let $p\geq 0$ with $p\geq q$, $r\geq 1$ and $\varepsilon>0$. Set $m=2p-2q\geq 0$. By the claim there are $l\geq 1$ and $k\geq l$ such that \eqref{eq-xk1} holds. We set
\[
n=p-q+n_k-n_l\geq 0.
\]
Then 
\begin{equation}\label{eq-p-n}
\|e_{p-n}\|_r=\|e_{q+n_l-n_k}\|_r<\varepsilon
\end{equation}
and
\begin{equation}\label{eq-p+n}
\|e_{p+n}\|_r = \|e_{q+n_k-n_l+2p-2q}\|_r= \|e_{q+n_k-n_l+m}\|_r < \varepsilon.
\end{equation}
Altogether, for any $p\geq 0$ with $p\geq q$, any $r\geq 1$ and $\varepsilon>0$, we have found some $n\geq 0$ such that \eqref{eq-p-n} and \eqref{eq-p+n} hold. By Theorem \ref{t-hcbil}, $B$ is hypercyclic; note that it is clearly enough by the continuity of $B$ that the condition holds for all large $p$.
\end{proof}

Since weakly bounded set are bounded, as are weak*-bounded sets by the Banach-Steinhaus theorem, we have the following in view of Example \ref{ex-bc}.

\begin{corollary}\label{c-weakseqlp}
Let $B_w$ be a weighted backward shift on $X=\ell^p(\mathbb Z)$, $1\leq p<\infty$. Then the following assertions are equivalent:
\begin{enumerate}
\item[\rm (i)] $B_w$ is hypercyclic;
\item[\rm (ii)] $B_w$ has an orbit with a non-zero weak sequential limit point in $X$;
\item[\rm (iii)] $B_w$ has an orbit with a non-zero weak$^*$-sequential limit point in $X$.
\end{enumerate}
\end{corollary}

Example \ref{ex-c0} shows that the corollary does not hold for $c_0(\Z)$. It does hold for the reflexive space $H(\C^*)$, which we already know by Theorem \ref{t-SchurMontel}. 

We finish our discussion of bilateral weighted backward shifts by showing that weak sequential hypercyclicity always implies hypercyclicity if $(e_n)_n$ is a basis. This result was obtained by B\`es, Chan and Sanders \cite{BCS05} for $\ell^p(\Z)$, $1\leq p<\infty$, and $c_0(\Z)$. However, their proof can easily be adapted to the general situation.  

\begin{theorem}\label{t-coobil}
Let $X$ be a Fréchet sequence space over $\Z$ in which $(e_n)_{n}$ is a basis, and let $B_w$ be a weighted backward shift on $X$. 

Suppose that there is some $x\in X$ such that, for any $y\in X$ there is some $(n_k)_k$ such that
\[
\{B_w^{n_k}x:k\geq 1\} \text{ is bounded and } B_w^{n_k}x\to y \text{ in } \K^{\Z}.
\]
Then $B_w$ is hypercyclic.
\end{theorem}

\begin{proof} We can assume that the shift is unweighted.

For this proof it is again convenient to work with seminorms. Thus let $(\|\cdot\|_r)_r$ be an increasing sequence of seminorms that defines the topology of $X$.

We will show that the characterizing condition of Theorem \ref{t-hcbil} is satisfied, that is, that for any $p\geq 0$, $r\geq 1$ and $\varepsilon >0$ there is some $n\geq 0$ such that
\begin{equation}\label{eq-hc}
\|e_{p+n}\|_r<\varepsilon\ \ \text{and}\text\ \ \|e_{p-n}\|_r<\varepsilon.
\end{equation}

First, since the mappings $z\mapsto z_ne_n$, $n\in \Z$, are equicontinuous on $X$ by the basis assumption, there is a constant $C>0$ and some $s\geq 1$ such that, for any $z\in X$,
\begin{equation}\label{eq-basis}
\|z_ne_n\|_r\leq C\|z\|_s, \ \ n\in \Z.
\end{equation} 

Now let $x\in X$ satisfy the hypothesis of the theorem. Taking $y=e_p$ we see that there exists a sequence $(n_k)_k$ and a constant $K>0$ such that 
\begin{equation}\label{eq-coohc}
\|B^{n_k}x\|_s\leq K, k\geq 1, \ \text{ and }\  x_{p+n_k}= [B^{n_k}x]_p \to 1.
\end{equation}
Let $\rho>\frac{CK}{\varepsilon}$. Applying the hypothesis again, this time with $y=2\rho e_p$, we find some $m\geq 0$ such that
\[
|x_{p+m}-2\rho| = |[B^{m}x-2\rho e_p]_p|< \rho,
\]
so that $|x_{p+m}|>\rho$. Using \eqref{eq-basis} and \eqref{eq-coohc} we have that, whenever $n_k\geq m$,
\begin{align*}
\|e_{p-(n_k-m)}\|_r &<\frac{1}{\rho} \|x_{p+m}e_{p-n_k+m}\|_r\leq \frac{C}{\rho}\Big\|\sum_{n\in\Z} x_{n+m}e_{n-n_k+m}\Big\|_s\\
&= \frac{C}{\rho}\Big\|\sum_{n\in\Z} x_{n}e_{n-n_k}\Big\|_s =\frac{C}{\rho}\|B^{n_k}x\|_s <\varepsilon.
\end{align*}
On the other hand, since $x\in X$ and $(e_n)_n$ is a basis,
\[
x_{p+n_k}e_{p+n_k}\to 0.
\]
By continuity of $B$, and by \eqref{eq-coohc}, there is some $n_k\geq m$ such that
\[
\|x_{p+n_k}e_{p+n_k-m}\|_r = \|B^m(x_{p+n_k}e_{p+n_k})\|_r<\frac{\varepsilon}{2}
\]
and 
\[
|x_{p+n_k}-1| <\frac{1}{2},
\]
so that $|x_{p+n_k}|>\frac{1}{2}$, and hence
\[
\|e_{p+(n_k-m)}\|_r = \frac{1}{|x_{p+n_k}|}\|x_{p+n_k}e_{p+n_k-m}\|_r < \varepsilon.
\]

Altogether, \eqref{eq-hc} holds with $n=n_k-m \geq 0$. 
\end{proof}

We note that the proof does not use the full strength of the hypothesis: we have taken for $y$ only multiples of unit sequences $e_n$.

\begin{corollary}\label{c-weakseqhc}
Let $X$ be a Fr\'echet sequence space over $\Z$ in which $(e_n)_n$ is a basis, and let $B_w$ be a weighted backward shift on $X$. 

\emph{(a)}Then the following assertions are equivalent:
\begin{enumerate}
\item[\rm (i)] $B_w$ is hypercyclic;
\item[\rm (ii)] $B_w$ is weakly sequentially hypercyclic.
\end{enumerate}

\emph{(b)} If, in addition, 
$X=Z^*$ is the dual of a Banach sequence space $Z$ under the natural pairing, then the following assertions are equivalent:
\begin{enumerate}
\item[\rm (i)] $B_w$ is hypercyclic;
\item[\rm (ii)] $B_w$ is weak*-sequentially hypercyclic.
\end{enumerate}
\end{corollary}

Let us recall that the notions of weak and weak*-sequential hypercyclicity are taken in the sense explained in the Introduction (see also Remark \ref{rem-seqhc}(b)).

\section{Zero-one laws for adjoint multiplication operators}\label{s-AdjMult}

One might wonder whether there are other natural classes of operators for which a zero-one law of orbital limit points holds. In this spirit, Chan and Seceleanu \cite{ChSe10} have considered adjoints of multiplication operators on Hilbert spaces of analytic functions. Godefroy and Shapiro \cite{GoSh91} had previously shown that they provide a rich source of hypercyclic operators. 

Chan and Seceleanu were able to obtain a zero-one law for the adjoints of multiplication operators (necessarily with non-constant multiplier) if the underlying space is a Bergman space $A^2(\Omega)$ on an arbitrary non-empty connected open set $\Omega$ in $\C$. But they note that their argument breaks down for the Hardy space $H^2(\D)$.

On the other hand, for the spaces considered by Godefroy and Shapiro, including $A^2(\Omega)$ and $H^2(\D)$, Costakis and Parissis \cite[Proposition 6.1]{CoPa12} have shown that, for adjoints of multiplication operators with non-constant multipliers, recurrence implies hypercyclicity. This suggests that a zero-one law might also hold on the Hardy space.

The aim of this section is not only to confirm this, but to obtain a zero-one law for the adjoints of multiplication operators on all the Hilbert spaces of analytic functions of one variable considered by Godefroy and Shapiro. Thanks to a recent result of Rong \cite{Ron22+} we are able to widen the framework to reflexive Banach spaces of analytic functions. 

Thus let $\Omega\subset \C$ be a non-empty connected open set and $X\neq \{0\}$ be a reflexive Banach space of analytic functions on $\Omega$.
In the sequel we will assume that, for any $z\in\Omega$, the point evaluation $X\to\C$, $f\mapsto f(z)$, is continuous and that every bounded analytic function $\phi$ on $\Omega$ defines a multiplication operator 
\[
M_\phi: X\to X, f\mapsto \phi f,
\]
such that $\|M_\phi\|\leq \sup_{z\in\Omega}|\phi(z)|$. Note that, by the closed graph theorem, any multiplication operator is continuous. 

Rong \cite{Ron22+} (but see also \cite[Exercise 4.4.5]{GrPe11}) has shown that if $\phi$ is non-constant then the adjoint operator, briefly called adjoint multiplier, 
\[
M^*_\phi : X^*\to X^*
\]
is hypercyclic if and only if
\[
\phi(\Omega)\cap\T\neq\varnothing,
\]
where $\T$ denotes the unit circle. Strictly speaking one has to distinguish between the Hilbert space adjoint and the Banach space adjoint, but they have the same dynamics, see \cite[Remark A.7]{GrPe11}.

Let us recall the definitions of the spaces of greatest interest to us. The Bergman space $A^2(\Omega)$ consists of all analytic functions $f$ on $\Omega$ such that $\|f\|^2:=\int_\Omega |f|^2 dA<\infty$, where $dA$ denotes area measure. The Hardy spaces $H^p(\D)$, $1\leq p<\infty$, consist of all analytic functions $f$ on the unit disk $\D$ such that $\|f\|^p:=\sup_{0\leq r< 1}\frac{1}{2\pi}\int_0^{2\pi} |f(re^{i\theta})|^p d\theta<\infty$; they are reflexive for $p>1$.

We have the following zero-one law.

\begin{theorem}\label{t-01adjoint}
Let $\Omega\subset \C$ be a non-empty connected open set and $X\neq \{0\}$ a reflexive Banach space of analytic functions that satisfies the assumptions stated above. Let $\phi$ be a non-constant bounded analytic function on $\Omega$, and let $M^*_{\phi}$ be the corresponding adjoint multiplier on $X^*$. Then the following assertions are equivalent:
\begin{enumerate}
\item[\rm (i)] $M^*_{\phi}$ is hypercyclic;
\item[\rm (ii)] $M^*_{\phi}$ has an orbit with a non-zero limit point;
\item[\rm (iii)] there exists a non-zero vectors $y\in X^*$ such that, for every neighbourhood $V$ of $y$, there is some $x\in X^*$ such that $(M^*_{\phi})^nx\in V$ for infinitely many $n\geq 0$;
\item[\rm (iv)] there are non-zero vectors $x,y\in X^*$ and some $0<r<\min(\|x\|,\|y\|)$ such that $(M^*_{\phi})^nx\in D(y,r)$ for infinitely many $n\geq 0$.
\end{enumerate}
\end{theorem}

This includes the zero-one law by Chan and Seceleanu for the Bergman space, and provides the missing one for the Hardy space $H^2(\D)$, and indeed for any $H^p(\D)$, $1<p<\infty$. 

\begin{proof}
It suffices to prove that (iv) implies (i). By way of contradiction, let us suppose that $M^*_{\phi}$ is not hypercyclic. Then, by the result mentioned above, $\phi(\Omega)\cap\T\neq\varnothing$. Since $\phi(\Omega)$ is connected, it must lie entirely inside the unit disk $\D$, or entirely outside $\overline{\D}$. 

Moreover, there exists some $\lambda\in \C$ such that $\lambda \phi(\Omega)$ intersects the unit circle. Again by the cited result, 
$\lambda M^{*}_{\phi} = M^{*}_{\lambda \phi}$ is hypercyclic, so that $M^{*}_{\phi}$ is supercyclic; see \cite{BaMa09} and \cite{GrPe11} for this notion. 

We now distinguish the two mentioned cases. Firstly, suppose that $\phi(\Omega)\subset \D$. Then
\[
\|M^{*}_{\phi}\| = \|M_{\phi}\|\le \sup_{z\in \Omega}|\phi(z)| \le 1,
\]
and hence $M^*_{\phi}$ is power bounded, which means that $(\|(M^*_\phi)^n\|)_n$ is bounded. Since $M^*_{\phi}$ is supercyclic and power bounded, a result by Ansari and Bourdon \cite{AnBo97} shows that all orbits under $M_\phi^*$ tend to zero. In particular, no orbit of $M^*_{\phi}$ can visit infinitely often a set $D(y,r)$ with $r<\|y\|$.

Secondly, suppose that $\phi(\Omega)\subset  \C\setminus \overline{\D}$. Then the function $\psi :=\frac{1}{\phi}$ is a bounded analytic function on $\Omega $ with $\psi(\Omega)\subset \D$. As above this implies that $\|M^*_{\psi}\|\leq 1$, and also that $M^*_{\psi}$ is supercyclic and power bounded, so that again all orbits of $M^*_{\psi}$ tend to zero.

By assumption there are $x,y\in X^*$ and $0<r\leq\|x\|$ such that $(M^*_{\phi})^nx\in D(y,r)$ for infinitely many $n\in\N_0$. Since $M_\psi^* M_\phi^*$ is the identity, we have that for these $n$,
\[
\|(M^*_{\psi})^ny-x\|\le \|(M^*_{\psi})^n\|\cdot\|y-(M^*_{\phi})^nx\|<r, 
\]
which is a contradiction because $(M^*_{\psi})^ny\to 0$ as $n\to\infty$ and $r\leq\|x\|$.
\end{proof}

\begin{corollary}\label{c-perrecajd}
Let $\Omega\subset \C$ be a non-empty connected open set and $X\neq \{0\}$ a reflexive Banach space of analytic functions that satisfies the assumptions stated above. Let $\phi$ be a non-constant bounded analytic function on $\Omega$, and let $M^*_{\phi}$ be the corresponding adjoint multiplier on $X^*$. 

\emph{(a)} If $M_\phi^*$ has a non-trivial periodic point, then it is hypercyclic.

\emph{(b)} The following assertions are equivalent:
\begin{enumerate}
\item[\rm (i)] $M_\phi^*$ is hypercyclic;
\item[\rm (ii)] $M_\phi^*$ is recurrent;
\item[\rm (iii)] $M_\phi^*$  has a non-zero recurrent vector.
\end{enumerate}
\end{corollary}

Part (b) improves the recurrent part of \cite[Proposition 6.1]{CoPa12}.

For adjoint multipliers on the Bergman space $A^2(\Omega)$ something more can be said than Theorem \ref{t-01adjoint}. We obtain here a trichotomy in the spirit of Theorem \ref{t-0-1-unilat_hyp}. Note that $A^2(\Omega)^*$ can be identified with $A^2(\Omega)$.

\begin{theorem}\label{teo: trichotomy bergman}
Let $\phi$ be a non-constant bounded analytic function on a non-empty connected open subset $\Omega\subset \mathbb C$, and let $M^*_{\phi}$ be the corresponding adjoint multiplier on $A^2(\Omega)$. Then the following trichotomy holds:
\begin{tabbing}
either  \=\kill
either  \> -- $M^*_{\phi}$ is hypercyclic,\\
or      \> -- for any $f\in A^2(\Omega)$, $(M^*_{\phi})^nf\to 0$,\\
or      \> -- for any non-zero $f\in A^2(\Omega)$, $\|(M^*_{\phi})^nf\|\to \infty$.
\end{tabbing}
\end{theorem}

\begin{proof}
Let us suppose that $M^*_{\phi}$ is not hypercyclic. Then we have seen in the previous proof that either $\phi(\Omega)\subset \D$ or $\phi(\Omega)\subset \C\setminus \overline{\D}$, and that in the former case, for any $f\in A^2(\Omega)$, $(M^*_{\phi})^nf\to 0$.

Now if $\phi(\Omega)\subset \C\setminus \overline{\D}$, then $\psi =\frac{1}{\phi}$ is a bounded analytic function on $\Omega $, so that $M_\psi$ defines an operator on $A^2(\Omega)$. Since $|\psi(z)|<1$ for all $z\in \Omega$, we have by the dominated convergence theorem, for any $g\in A^2(\Omega)$, 
\[
\|M^n_{\psi}g\|^2=\int_\Omega |\psi(z)^ng(z)|^2 dA\to 0.
\]

Let $f\in A^2(\Omega)^*=A^2(\Omega)$ be any non-zero function. Take $g\in A^2(\Omega)$ such that $\langle g, f\rangle\ne 0$. Since $M_\psi g\neq 0$ and $M_\phi M_\psi$ is the identity, we have that
\[
\Big\langle \frac{M^n_{\psi}g}{\|M^n_{\psi}g\|},(M^*_{\phi})^nf\Big\rangle =\frac{1}{\|M^n_{\psi}g\|}\langle g,f\rangle  \to \infty.
\]
This implies that $\|(M^*_{\phi})^nf\|\to \infty$.
\end{proof}

As a consequence, we obtain for the Bergman spaces the following improvement of the zero-one law for adjoint multipliers. Recall that $x\in X$ is an irregular vector for an operator $T$ on $X$ if $\liminf_{n\to\infty} \|T^nx\|=0$ and $\limsup_{n\to\infty} \|T^nx\|=\infty$.

\begin{corollary}\label{coro:0-1 bergman}
Let $\phi$ be a non-constant bounded analytic function on a non-empty connected open subset $\Omega\subset \mathbb C$, and let $M^*_{\phi}$ be the corresponding adjoint multiplier on $A^2(\Omega)$. Then the following assertions are equivalent:
\begin{enumerate}
\item[\rm (i)] $M^*_{\phi}$ is hypercyclic;
\item[\rm (ii)] $M^*_{\phi}$ has an orbit with a non-zero limit point;
\item[\rm (iii)] there is some $f\in A^2(\Omega)$ whose orbit has infinitely many members in an open ball whose closure avoids the origin;
\item[\rm (iv)] $M^*_{\phi}$ has an irregular vector.
\end{enumerate}
\end{corollary}

The characterizing condition (iii) was already obtained by Chan and Seceleanu \cite[Theorem 2.1]{ChSe10} and strengthens condition (ii) in Theorem \ref{t-01adjoint}. We do not know if the same improvement is possible on general spaces.

\begin{question} 
In the setting of Theorem \ref{t-01adjoint}, is $M^*_{\phi}$ hypercyclic as soon as there is some $x\in X$ whose orbit has infinitely many members in an open ball whose closure avoids the origin?
\end{question} 

On the other hand, the equivalence of (i) and (iv) was already obtained by Bernardes et al. \cite[Theorems 9 and 26]{BBMP15}; their proof even extends to all spaces in the setting of Theorem  \ref{t-01adjoint}.

\begin{remark}
Incidentally, the existence of an irregular vector for a unilateral weighted backward shift on $\ell^p$, $1\leq p<\infty$, or $c_0$ does not imply hypercyclicity, see \cite[Theorem 9, Remark 24]{BBMP15}.
\end{remark}

For the Hardy spaces we can also prove a trichotomy if we assume that the symbol of the adjoint multiplier is analytic on a neighbourhood of the disk. Note that if $M^*_\phi$ is defined on $H^p(\D)$, then the underlying multiplication operator $M_\phi$ acts on $H^{q}(\D)$, where $q$ is the conjugate exponent.

\begin{theorem}\label{teo: trichotomy hardy}
Let $\phi$ be a non-constant analytic function on a neighbourhood of $\overline{\D}$, and let $M^*_{\phi}$ be the corresponding adjoint multiplier on $H^p(\D)$, $1<p<\infty$. Then the following trichotomy holds:
\begin{tabbing}
either  \=\kill
either  \> -- $M^*_{\phi}$ is hypercyclic,\\
or      \> -- for any $f\in H^p(\mathbb D)$, $(M^*_{\phi})^nf\to 0$,\\
or      \> -- for any non-zero $f\in H^p(\mathbb D)$, $\|(M^*_{\phi})^nf\|\to \infty$.
\end{tabbing}
\end{theorem}

\begin{proof}
If $M^*_{\phi}$ is not hypercyclic then  either $\phi(\D)\subset \D$ or $\phi(\D)\subset \C\setminus \overline{\D}$.
Assume first that $\phi(\D)\subset {\D}$. Note that, since $H^q(\D)^*$ and $H^p(\D)$ are isomorphic but not isometric, we may have $\|M_\phi^*\|>1$ when considered as an operator acting on $H^p(\D)$. We will thus consider $M_\phi^*$  as an operator on $H^q(\D)^*$, where it is a power bounded operator. Then we can apply the Ansari and Bourdon result from \cite{AnBo97} as we did in Theorem \ref{t-01adjoint} to conclude that, for any $f\in H^q(\mathbb D)^*$, $(M^*_{\phi})^nf\to 0$ or equivalently, $(M^*_{\phi})^nf\to 0$ for any $f\in H^p(\mathbb D)$.

Let us now show that if $\phi(\D)\subset \C\setminus \overline{\D}$ then the third alternative holds. Thus assume that $|\phi(z)|>1$ for all $z\in\mathbb D$. 

We first note that $\phi(\T)\cap \mathbb T$ is a finite set. Indeed, the function $t\mapsto h(t):= |\phi(e^{it})|^2=\phi(e^{it})\overline{\phi(e^{it})}$ is a real analytic function of $\cos(t)$ and $\sin(t)$, so real analytic on $\R$. If $\phi(\T)\cap \mathbb T$ were not finite, then $h$ would be constant $1$ by the identity principle. But then $\phi(\T)\subset\T$, which contradicts the maximum principle since $\phi(\D)\subset \C\setminus \overline{\D}$. 

Let now $\psi=\frac{1}{\phi},$ which is bounded and analytic on $\D$, so that $M_\psi$ defines an operator on $H^{q}(\mathbb D)$. Since any $g\in H^{q}(\mathbb D)$ extends to a function in $L^q(\T)$ through radial limits, which exist almost everywhere, we have that
\[
\|M_\psi^ng\|^q=\frac{1}{2\pi} \int_0^{2\pi} |\psi(e^{i\theta})|^{nq}|g(e^{i\theta})|^{q}d\theta.
\]
Since, by the previous argument, $|\psi|<1$ on $\mathbb T$ with at most finitely many exceptions, the dominated convergence theorem implies that $M_\psi^ng\to 0$. As in the proof of Theorem \ref{teo: trichotomy bergman}, we then deduce $\|(M_\phi^{*})^nf\| \to \infty$ for any non-zero $f\in H^q(\mathbb D)^*$, or equivalently $\|(M_\phi^{*})^nf\| \to \infty$ for any non-zero $f\in H^p(\mathbb D)$.
\end{proof}

We mention that, as a consequence of Theorem \ref{teo: trichotomy hardy}, Corollary \ref{coro:0-1 bergman} also holds on Hardy spaces if $\phi$ is a non-constant analytic function on a neighbourhood of $\overline{\D}$.

\begin{question}
(a) Let $\phi(z)=e^{\pi/2}\left(\frac{1+z}{1-z}\right)^i$. Then $\phi(\mathbb D)$ is the annulus $\{1<|z|<e^{\pi}\}$. Does the orbit of every non-zero function $f\in H^2(\mathbb D)$ under $M_\phi^*$ go to infinity? We mention that $M_\phi^*$ is weakly hypercyclic by an unpublished result of Shkarin \cite[Theorem 1.9]{Shk12}.

(b) Does Theorem \ref{teo: trichotomy hardy} hold for any $\phi\in A(\D)$, that is, if $\phi$ is continuous on $\overline{\D}$ and analytic in $\D$?
\end{question}

There is a close link between adjoint multipliers on $H^2(\mathbb D)$ and functions $\phi(B)$ of the (unweighted) backward shift $B$ on $\ell^2=\ell^2(\N_0)$ if $\phi$ is analytic on a neighbourhood of $\overline{\D}$, see \cite[Proposition 4.41]{GrPe11}. It is therefore interesting to study zero-one laws for $\phi(B)$ on any of the complex spaces $\ell^p$, $1\leq p <\infty$, and $c_0$. 

\begin{corollary}
Let $\phi$ be a non-constant analytic function on a neighbourhood of $\overline{\mathbb D}$. Let $B$ be the backward shift on a complex sequence space $X=\ell^p$, $1\leq p \leq 2$. Then the following assertions are equivalent:
\begin{enumerate}
\item[\rm (i)] $\phi(B)$ is hypercyclic;
\item[\rm (ii)] $\phi(B)$ has an orbit with a non-zero limit point;
\item[\rm (iii)] $\phi(\D)\cap \T\neq\varnothing$.
\end{enumerate}
\end{corollary}

\begin{proof}
Since (iii) implies (i) by \cite[Theorem 4.43]{GrPe11}, it suffices to show that (ii) implies (iii). Suppose that $\phi(B)$, as an operator on $\ell^p$, has an orbit with a non-zero limit point. By continuous inclusion of $\ell^p$ in $\ell^2$, this is also true for $\phi(B)$ as an operator on $\ell^2$. Since, by \cite[Proposition 4.41]{GrPe11}, $\phi(B)$ on $\ell^2$ is conjugate to the (Banach space) adjoint $M_\phi^*$ on $H^2(\mathbb D)$, $M_\phi^*$ has an orbit with a non-zero limit point in $H^2(\mathbb D)$. An application of Theorem \ref{teo: trichotomy hardy} shows that $M_\phi^*$ is hypercyclic on $H^2(\mathbb D)$, hence that $\phi(\D)\cap\T\neq\varnothing$, which had to be shown.
\end{proof}

We do not know if the result also holds for $p>2$ or on $c_0$.

\begin{question}\label{q-phiB}
Let $\phi$ be a non-constant analytic function on a neighbourhood of $\overline{\D}$. Consider $\phi(B)$ on a complex sequence space $X=\ell^p$, $2< p <\infty$, or on $c_0$. Is $\phi(B)$ hypercyclic as soon as it has a non-zero limit point? 
\end{question}

\end{document}